\documentclass{amsart}

\usepackage{amssymb}
\usepackage[all]{xy}

\newtheorem{thm}{Theorem}

\newtheorem{lem}[thm]{Lemma}
\newtheorem{cor}[thm]{Corollary}

\newtheorem{prop}[thm]{Proposition}

\newtheorem{complem}[thm]{Complement}
   
\theoremstyle{definition}
\newtheorem{defn}[thm]{Definition}

\newtheorem{say}[thm]{}
\newtheorem{exmp}[thm]{Example}

\newtheorem{prob}[thm]{Problem}
\newtheorem{ques}[thm]{Question}    

\newtheorem{rem}[thm]{Remark}          

\newtheorem*{ack}{Acknowledgments}      
\newtheorem{notation}[thm]{Notation}   
  
\newtheorem{defn-thm}[thm]{Definition--Theorem}  
\newtheorem{defn-lem}[thm]{Definition--Lemma}  

\theoremstyle{remark}

\renewcommand{\c}[0]{{\mathbb C}}  

\renewcommand{\o}[0]{{\mathcal O}} 
\newcommand{\z}[0]{{\mathbb Z}}

\renewcommand{\r}[0]{{\mathbb R}} 

\renewcommand{\a}[0]{{\mathbb A}}

\newcommand{\s}[0]{{\mathbb S}}

\newcommand{\p}[0]{{\mathbb P}}
\newcommand{\f}[0]{{\mathbb F}}
\newcommand{\q}[0]{{\mathbb Q}}
\newcommand{\map}[0]{\dasharrow}
\newcommand{\qtq}[1]{\quad\mbox{#1}\quad}
\newcommand{\spec}[0]{\operatorname{Spec}}

\newcommand{\rank}[0]{\operatorname{rank}}

\newcommand{\aut}[0]{\operatorname{Aut}}

\newcommand{\disc}[0]{\operatorname{discr}} 
 
\newcommand{\chow}[0]{\operatorname{Chow}}
\newcommand{\chr}[0]{\operatorname{char}}

\newcommand{\hilb}[0]{\operatorname{Hilb}}
\newcommand{\grass}[0]{\operatorname{Grass}}

\newcommand{\univ}[0]{\operatorname{Univ}}

\newcommand{\onto}[0]{\twoheadrightarrow}

\newcommand{\tr}[0]{\operatorname{tr}}

\newcommand{\Div}[0]{\operatorname{Div}}

\newcommand{\rc}[0]{\operatorname{RatCurve}}

\def\into{\DOTSB\lhook\joinrel\to}

\def\loccoh#1.#2.#3.#4.{H^{#1}_{#2}(#3,#4)}

\DeclareMathAlphabet{\mathchanc}{OT1}{pzc}%
                                {m}{it}

\newcommand{\simb}[0]{\stackrel{bir}{\sim}}
\newcommand{\sims}[0]{\stackrel{stab}{\sim}}

\newcommand{\PGL}{\mathrm{PGL}}

\newcommand{\SO}{\mathrm{SO}}
\newcommand{\OO}{\mathrm{O}}
\newcommand{\PGO}{\mathrm{PGO}}

\newcommand{\OG}{\mathrm{OG}}
\newcommand{\OF}{\mathrm{OF}}

\newcommand{\PTC}{\mathrm{PTC}}

\newcommand{\sym}[0]{\operatorname{Sym}}

\usepackage[all]{xy}\xyoption{dvips}

\newcommand{\norm}[0]{\operatorname{norm}}

\newcommand{\maps}[0]{\operatorname{Map}}

\newcommand{\mo}[0]{\operatorname{M}^{\rm \circ}}

\begin{document}
\bibliographystyle{amsalpha}


 \title{Quadratic solutions of quadratic forms}
 \author{J\'anos Koll\'ar}

 \maketitle

The aim of this note is to
study solutions of a homogeneous
quadratic equation $q(x_0,\dots, x_{n+1})=0$, defined over a field $k$,  where the
$x_i=h_i(u_0,\dots, u_r)$ are themselves  homogeneous polynomials of a given degree $d$.  
Equivalently, we are looking at rational maps  defined over $k$, 
$$
\Phi: \p^r\map Q^n:=\bigl(q(x_0,\dots, x_{n+1})=0\bigr)\subset \p^{n+1}
$$
from projective $r$-space to an $n$-dimensional quadric hypersurface.  
We denote the set of such maps by
$\maps_d(\p^r, Q^n)(k)$. It is easy to see that there is a
projective algebraic set $\maps_d(\p^r, Q^n)$, defined over $k$, whose 
$K$ points can be naturally identified with $\maps_d(\p^r, Q^n)(K)$
for any field extension $K\supset k$.  Thus we aim to describe
the algebraic sets $\maps_d(\p^r, Q^n)$ and give concrete parametrization of their $k$ points.

If $d=1$ then the image of $\Phi$ is a linear subspace of $Q^n$, hence
the study of linear solutions of quadratic equations is essentially equivalent
to the theory of   orthogonal Grassmannians. These are homogeneous spaces under the orthogonal group  $\aut(Q^n)=\PGO(q)$ and quite well understood,  see 
Section \ref{sec.OG}.

Another extreme case is when $r=1$. Solving 
$q(x_0,\dots, x_{n+1})=0$,  where the
$x_i=h_i(u,v)$ are   homogeneous polynomials in  2 variables $u,v$
is  equivalent to  studying maps $\p^1\to Q^n$.
In algebraic geometry there has been considerable interest in understanding
 maps from $\p^1$ to a given variety. Most of these works deal with
more general varieties over algebraically closed fields. Quadrics are special cases of  homogeneous spaces \cite{MR1882330}
and of low degree hypersurfaces \cite{MR2070144}. 
See the introductory notes   \cite{ar-ko, k-lisbon} or the more complete 
 treatment  given in \cite{rc-book} for details. 
The existence of  maps  $\p^1\to X$ over finite fields is studied in \cite{k-looking},
but not much is known  about the structure of these spaces over arbitrary fields. For quadrics, we prove the following in Section \ref{sec.2}.

\begin{thm} \label{curve.map.thm}  
Let $Q^n$ be a smooth quadric of dimension $n\geq 3$. Then
$$
\maps_{d}(\p^1, Q^n)\simb 
\left\{
\begin{array}{l}
Q^n\times \p^{nd}\qtq{if $d$ is even, and}\\
\OG(\p^1, Q^n)\times \p^{nd-n+3}\qtq{if $d$ is odd,}
\end{array}
\right.
$$
where 
$\simb$ denotes birational equivalence and 
$\OG(\p^1,Q^n)$ the  orthogonal Grassmannian of lines in $Q^n$.
\end{thm}

Instead of parametrizing maps $\p^1\to Q^n$, it is also of interest to
parametrize spaces of rational curves contained in $Q^n$.
For odd degrees, the latter problem is almost equivalent to
Theorem \ref{curve.map.thm}; see
Corollary \ref{RC.quad.odd.thm}. However, for even degrees
parametrizing spaces of  rational curves
 is much harder and we get full answers in only a few cases; see 
  Theorem \ref{RC.quad.thm} and Proposition \ref{06.16.prop}.

The rest of the paper is devoted to the next case  $r=d=2$, which  is considerably harder. Thus we aim to understand  rational maps   
$$
\Phi: \p^2\map Q^n:=\bigl(q(x_0,\dots, x_{n+1})=0\bigr)\subset \p^{n+1}
$$
defined by  degree $2$ homogeneous polynomials  $x_i=h_i(u,v, w)\in k[u,v,w]$.
As in Theorem   \ref{curve.map.thm}, our aim is to describe 
the moduli spaces $\maps_2(\p^2, Q^n) $ up to birational equivalence, using
$Q^n$ and its orthogonal Grassmannians.
If $Q^n$ is singular then, after a suitable linear coordinate change, we can eliminate some of the variables $x_i$ from $q$. Thus from now on $Q^n$ denotes a smooth quadric hypersurface of dimension $n$ over a field $k$.

With the exception of Sections \ref{sec.2}--\ref{sec.3} we  assume form now on that  $\chr k\neq 2$, though this is not always necessary.

In low dimensions we prove in  (\ref{conics.say}.2), (\ref{q.2-q20.say}.2) and  (\ref{3f.say}.3)
 that
$$
\begin{array}{rcl}
\maps_2(\p^2, Q^1)&\simb &Q^1 \times\p^4,\\
\maps_2(\p^2, Q^2)&\simb & Q^2\times \p^{8},\\
\maps_2(\p^2, Q^3)&\simb &\bigl(Q^3\amalg  \OG(\p^1, Q^3)\bigr)\times \p^{11}.
\end{array}
$$

Let $\maps_2^{\circ}(\p^2, Q^n)\subset \maps_2(\p^2, Q^n)$ denote the
open subset consisting of maps that are everywhere defined.
If $\Phi\in \maps_2^{\circ}(\p^2, Q^n)$ and we work over $\c$
then the image of the fundamental class
$\Phi_*[\p^2]$ is in $H_4(Q^n, \z)$. 
If $n= 2$ (resp.\ $3$)  then  $H_4(Q^n,\z)\cong \z$ is generated by
$Q^{2}$ (resp.\ $Q^3\cap H$) where $H$ denotes the hyperplane class.
 If $n\geq  5$   then  $H_4(Q^n,\z)\cong \z$ but this time it is generated by
the class of any 2-plane  $L^2\subset Q^n$. 
 The most interesting case is  $n=4$. Then 
 $H_4(Q^4,\z)\cong \z[A]+\z[B]$ where $A, B\subset Q^4$ are disjoint
planes; see \cite[Book IV, Sec.XIII.4]{hodge-ped} on linear spaces
in quadrics. The same holds over any algebraically closed field
if we replace  $H_4(Q^n,\z)$ by the group   $N_2(Q^n,\z)$ of 2-cycles modulo numerical equivalence. 
This suggests that the really interesting case is $n=4$.  

The key player in our study is the Veronese surface.

\begin{defn} A {\it Veronese surface} is the image of the map
$$
\p^2\to \p^5\qtq{given by}
(u:v:w)\mapsto (u^2:v^2:w^2:uv:vw: uw),
$$
composed with any automorphism of $\p^5$. 
They  were first studied  in \cite{veronese}. 
Veronese surfaces have numerous exceptional and extremal properties.

$\bullet$  If a surface contains a conic through any 2 of its points then it is 
a plane, quadric or Veronese; this follows from the Kronecker--Castelnuovo theorem; see  \cite{cast-KC}.

$\bullet$  If the  secant lines of a  surface  cover
 a subvariety of dimension $\leq 4$, then it is the Veronese
(or it is contained in a linear subspace of dimension $\leq 4$).
This was proved by Severi \cite{sev-veronese} 
(see \cite{MR825917} for a short proof) and 
generalized by Kuiper and Pohl \cite{MR0145545, MR0494122}
to non-algebraic  surfaces. 

$\bullet$  Veronese surfaces give rise to remarkable
birational maps. As pointed out by F.~Russo, 
this seems to have been known to Bertini
(see \cite[Secs.XV.9--10]{bertini-book} or
\cite[Secs.XVI.9--10]{bertini-book-g}) 
and Coble \cite[Thm.17]{MR1501210},
but the earliest explicit mention may be in 
 \cite[p.188]{MR0034048}. 
A modern treatment is given by  Ein and Shepherd-Barron \cite{ein-she}. 
These connections have been extended  by Pirio and Russo  to include  Jordan algebras \cite{MR3330541}. See the book \cite{rus-book} for a detailed treatment of many of these results.

We will use a classification of  projections and equations of Veronese surfaces; see
(\ref{veronese.pr.say}) and (\ref{vero.eq.say}).
\end{defn}

The geometric version of the main result is  the following;

\begin{thm}\label{p2.in.q4.C.thm}
 Let $Q^4$ be a smooth quadric 4-fold defined over an algebraically closed field.  Then $\maps_2(\p^2, Q^4)$ has 5 irreducible components, each birational to $\p^{20}$, and
$\maps_2^{\circ}(\p^2, Q^4)$ is dense in $\maps_2(\p^2, Q^4)$.
The 5 components of $\maps_2^{\circ}(\p^2, Q^4)$ can be described as follows
\begin{enumerate}
\item (Veronese)  In these cases $\Phi(\p^2)\subset Q^4$
is a Veronese surface. There are 2 such components corresponding to
 $\Phi_*[\p^2]=3[A]+[B]$ and $\Phi_*[\p^2]=[A]+3[B]$, where 
$A, B\subset Q^4$ are disjoint planes.
\item (Projected Veronese) In these cases $\Phi(\p^2)\subset Q^4$
is a singular projection of a Veronese surface. There is 1 such component and 
 $\Phi_*[\p^2]=2[A]+2[B]$.
\item (Quadruple plane)   In these cases $\Phi(\p^2)\subset Q^4$
is a plane.  There are 2 such components corresponding to
 $\Phi_*[\p^2]=4[A]$ and $\Phi_*[\p^2]=4[B]$.
\end{enumerate}
\end{thm}

The corresponding subvarieties of $\maps_2^{\circ}(\p^2, Q^4)$
will be denoted by
$\maps_{V}^{\circ}(\p^2, Q^4)$, $\maps_{PV}^{\circ}(\p^2, Q^4)$, 
$\maps_{QP}^{\circ}(\p^2, Q^4)$ and their closures in $\maps_2(\p^2, Q^4)$
by $\maps_{V}(\p^2, Q^4)$, $\maps_{PV}(\p^2, Q^4)$ and  
$\maps_{QP}(\p^2, Q^4)$. 

Our main aim is to describe $\maps(\p^2, Q^4) $ over  arbitrary fields $k$, so  Theorem \ref{p2.in.q4.C.thm} is but a  special case of the following 3 results when $k$ is algebraically closed. 
The most interesting is  the description of the Veronese locus; it is  proved in Section \ref{sec.veron}. 
The easier Propositions \ref{4-folds.2.thm}--\ref{4-folds.1.thm} are treated in Section \ref{sec.degmaps}.
See Notation \ref{basic.not.not} for the terminology involving quadratic forms.

\begin{thm}[Veronese]\label{4-folds.3.thm}   Let  $Q^4$  be a smooth quadric 4-fold over a field $k$. Then 
$$
\maps_{V}(\p^2, Q^4)\simb \bigl\{\pm\sqrt{-\Delta}\bigr\}\times Q^4\times \p^{16},
$$
where   $\Delta $ is the discriminant of $Q^4$ and  
$\bigl\{\pm\sqrt{-\Delta}\bigr\} $ denotes the 2-point algebraic set defined by
the equation $t^2+\Delta=0$ over $k$. 

Furthermore, $\maps_{V}^{\circ}(\p^2, Q^n)(k)$
is nonempty  iff   $-\Delta$ is a square in $k$ and $Q^4(k)\neq \emptyset$. If these hold  then
$\maps_{V}(\p^2, Q^4)$ is birational to $\p^{20}\amalg \p^{20}$, the disjoint union of 2 copies of $\p^{20}$
\end{thm}

The description of the projected Veronese maps is not hard and 
the quadruple planes are obvious. The following are 
 special cases of (\ref{pvs.map.say}.2) and 
(\ref{quadr.pl.say}.1).

\begin{prop}[Projected Veronese]\label{4-folds.2.thm} Let  $Q^4$  be a smooth quadric 4-fold.
Then
$$
\maps_{PV}(\p^2, Q^4)\simb  \OG(\p^1,Q^4)\times \p^{15},
$$
where $\OG(\p^1,Q^4)$ denotes the  orthogonal Grassmannian of lines in $Q^4$. 
 Furthermore, $\maps_{PV}^{\circ}(\p^2, Q^4)(k)$
is nonempty  iff the anisotropic rank of $Q^4$ is 
$\leq 2$ and then $\maps_{PV}(\p^2, Q^4)$ is birational to $\p^{20}$.
\end{prop}

\begin{prop}[Quadruple plane]\label{4-folds.1.thm}
 Let  $Q^4$  be a smooth quadric 4-fold. Then 
$$
\maps_{QP}(\p^2, Q^4)\simb \OG(\p^2,Q^4)\times \p^{17},
$$
where $\OG(\p^2,Q^4)$ denotes the  orthogonal Grassmannian of $2$-planes in $Q^4$.  Furthermore, $\maps_{QP}^{\circ}(\p^2, Q^n)(k)\neq \emptyset$
iff $Q^4$ is split over $k$  and then
$\maps_{QP}(\p^2, Q^4)$ is birational to $\p^{20}\amalg \p^{20}$.
\end{prop}

Note that if the anisotropic rank of $Q^4$ is $\leq 2$ and $-\Delta$ is a square 
then $Q^4$ is split. Thus for non-split forms we have either only
1  irreducible component with smooth $k$-points
(if the anisotropic rank is $2$) or 2 irreducible components with smooth $k$-points (if $-\Delta$ is a square).

The above results count maps from a fixed $\p^2$ to a  quadric
$Q^4$. 
Frequently it is more interesting to count maps 
up to coordinate changes by $\aut(\p^2)$ or by  $\aut(Q^4)$.

Taking quotient by  $\aut(\p^2)$ is essentially the same as
working with the components of the Chow variety $\chow(Q^4)$
that parametrize Veronese  surfaces  (resp.\ projected Veronese  surfaces).
Let us denote these by 
$\chow_{V}(Q^4)$  (resp.\  $\chow_{PV}(Q^4)$). 
Let  $ H_5$ denote the hyperplane class on $\p^5$ and 
let $F_k\subset Q^4_k$ be a surface that is a  Veronese  surface over $\bar k$.
Then $|2H_5|_{F_k}+K_{F_k}|$ is the linear system of lines giving an isomorphism $F_k\cong \p^2_k$. Therefore
 $\chow_{V}(Q^4)\times \aut(\p^2)\simb \maps_{V}(\p^2, Q^4)$ and
similarly for $\chow_{PV}(Q^4)$.
Hence the previous results determine the stable birational types
of $\chow_{V}(Q^4)$  and  $\chow_{PV}(Q^4)$. 
The proofs  yield the following slightly stronger forms;
see (\ref{Ver.final.pf.say}.6) and (\ref{pvs.map.say}.1).

\begin{cor}[Chow variety version]\label{4-folds.3.chow}   Let  $Q^4$  be a smooth quadric 4-fold over a field $k$. Then 
$$
\begin{array}{lcl}
\chow_{V}(Q^4)\times \p^4&\simb& \bigl\{\pm\sqrt{-\Delta}\bigr\}\times Q^4\times \p^{12}\qtq{and}\\
\chow_{PV}(Q^4)&\simb & \OG(\p^1,Q^4)\times \p^{7}.
\end{array}
$$
\end{cor}

It is quite likely that one can cancel the factor $\p^4$ but
my proof does not seem to work without it.

 $\chow_{PV}(Q^4) $ is contained in  the
irreducible component  $\chow_{CI}(Q^4) \subset \chow(Q^4) $ parametrizing complete intersections of $Q^4$ with a hyperplane and a quadric.
Another interesting irreducible component of $\chow(Q^4) $  parametrizing degree 4 surfaces is discussed in (\ref{other.comp.rem}). 

We show in Section \ref{sec.coo.ch} that working  
up to coordinate changes by both $\aut(\p^2)$ and  $\aut(Q^4)$ yields a very nice answer, especially in the Veronese case.

\begin{thm}[Veronese] \label{V.UP.TO.ISOM.THM}
Let $k$ be a field. A smooth quadric hypersurface $Q^4\subset \p^5$ contains a Veronese surface iff its equation can be written as
$$
Q^4_{a,b}=\bigl(x_0x_5=x_1^2+ax_2^2+bx_3^2+abx_4^4\bigr).
\eqno{(\ref{V.UP.TO.ISOM.THM}.1)}
$$
Up to coordinate changes by $\aut(\p^2)\times\aut\bigl(Q^4_{a,b}\bigr)$,
the quadric $Q^4_{a,b}$ contains a unique Veronese surface. 
A  representative  is given by
$$
(u{:}v{:}w)\mapsto 
\bigl(u^2+aw^2:
uv : vw : uw : w^2 : v^2+bw^2\bigr).
\eqno{(\ref{V.UP.TO.ISOM.THM}.2)}
$$
\end{thm}

\begin{complem} \label{V.UP.TO.ISOM.compl} Let $k$ be a field. Then
$Q^4_{a,b}\leftrightarrow  (au^2+bv^2+w^2=0)$ gives  a one-to-one correspondence between
\begin{itemize}
\item pairs $(V\subset Q^4)$---a Veronese surface contained in a smooth quadric hypersurface defined over $k$,  up to isomorphism---and
\item smooth plane conics  $C\subset \p^2$ defined over $k$, up to isomorphism. 
\end{itemize}
\end{complem}

A proof of  this is given in (\ref{vero.eq.say}.6).
Another  way of obtaining a conic  out of a pair  $(V\subset Q^4)$ is as follows. 
Pick $p\in V$. Intersecting $V$ with the tangent plane of $Q^4$ at $p$
gives a hyperplane section of $V$ that is  singular at $p$.
Since the hyperplane sections of $V$ are conics, we expect that
we get a pair of lines and there is a curve
$C_{a,b}\subset V$ where this intersection is a double line.

A direct computation using  (\ref{V.UP.TO.ISOM.THM}.2) 
(that I did not find illuminating) shows 
that 
$C_{a,b}=(bu^2+av^2+abw^2=0)$. We can rewrite this as
$b(au_1)^2+a(bv_1)^2+abw^2=ab(au_1^2+bv_1^2+w^2)$.

\begin{rem} \label{quat.form.rem}
We can fully diagonalize the quadric $Q^4_{a,b}$ by replacing
$x_0, x_5$ by $\tfrac12(x_0+x_5)$ and $\tfrac12(x_0-x_5)$.
Then we end up with  the identity
$$
\tfrac14\bigl(u^2+v^2+(a+b)w^2\bigr)^2-
\tfrac14\bigl(u^2-v^2+(a-b)w^2\bigr)^2
=
(uv)^2+a(vw)^2+b(uw)^2+ab(w^2)^2.
$$
Skopenkov pointed out that (if $a, b$ are real and positive) this is equivalent to the quaternionic identity
$$
\tfrac14 (|A|^2+|B|^2)^2-\tfrac14 (|A|^2-|B|^2)^2=|AB|^2
$$
for the quaternions
$A=u+iw\sqrt{a}$ and $B=v+jw\sqrt{b}$.
\end{rem}

The following 2 Propositions are also proved in Section \ref{sec.coo.ch}.

\begin{prop}[Projected Veronese]  \label{PV.UP.TO.ISOM.THM}
Let $k$ be a field. A smooth quadric hypersurface $Q^4\subset \p^5$ contains a projected Veronese surface iff its equation can be written as
$$
Q^4_{a}=\bigl(x_0x_1+x_2x_3=x_4^2+ax_5^2\bigr).
$$
After  coordinate changes by $\aut(\p^2)\times\aut\bigl(Q^4_a\bigr)$,
projected Veronese surfaces can be brought to the form
$$
(u{:}v{:}w)\mapsto \bigl(0: q(u{:}v{:}w): u^2+av^2: w^2: uw:vw\bigr)
$$
where $\deg q=2$.
Restricting $q$ to the line  $(w=0)$
 gives an isomorphism between the  moduli  of
 projected Veronese surfaces in $Q^4_a$
(up to  $\aut(\p^2)\times\aut(Q^4_a)$)
and  $\sym^2(\p^1)\big/\OO(u^2+av^2)$. 
\end{prop}

\begin{prop}[Quadruple plane]   \label{QP.UP.TO.ISOM.THM}
A smooth quadric hypersurface $Q^4\subset \p^5$ contains a 2-plane iff its equation can be written as
$$
Q^4_{\rm split}=\bigl(x_0x_3+x_1x_4+x_2x_5=0\bigr).
$$
If   $\chr k\neq 3$  then, after  coordinate changes by $\aut\bigl(Q^4_{\rm split}\bigr)$,
quadruple planes can be brought to the form
$$
(u{:}v{:}w)\mapsto 
\bigl(\partial C/\partial u: \partial C/\partial v: \partial C/\partial w:0:0:0\bigr),
$$
where $C=C(u,v,w)$ is a homogeneous cubic. 
This gives an isomorphism between the  moduli of quadruple planes
(up to  $\aut(\p^2)\times\aut\bigl(Q^4_{\rm split}\bigr)$)
and  the moduli  of
 plane cubic curves (up to $\aut(\p^2)$). 
\end{prop}

Next we move to quadrics of dimension $\geq 5$.
 Note that $\Phi(\p^2)$
is always contained in a 5-dimensional linear subspace of $\p^{n+1}$, so
the cases  $n\geq 5$ reduce to the $n\leq 4$ cases is principle
but  the  parametrization of the solutions for $n\geq 5$ 
 proceeds somewhat differently. The following is shown  in Section \ref{sec.5.dim}.

\begin{thm} \label{n>4.main.thm}
 Assume that $n\geq 5$ and   $Q^n(k)\neq \emptyset$. Then 
$$
\maps_2(\p^2, Q^n)\simb \sym^2\bigl(\OG(\p^2, Q^n)\bigr)\times \p^{14},
\eqno{(\ref{n>4.main.thm}.1)}
$$
where   $\sym^2 $ denotes the symmetric square.
\end{thm}

Using (\ref{pls.say}.2) we can rewrite (\ref{n>4.main.thm}.1) as
$$
\maps_2(\p^2, Q^n)\simb \sym^2\bigl(\OG(\p^1, Q^{n-2}_W)\bigr)\times \p^{2n+10},
\eqno{(\ref{n>4.main.thm}.2)}
$$
where  $Q^{n-2}_W$  is the Witt reduction of $Q^n$ (see Notation
\ref{basic.not.not}). 
This variant  actually shows the structure of $\maps_2(\p^2, Q^n) $  better than
(\ref{n>4.main.thm}.1). 

\begin{rem}
If $n=4$ then $ \bigl\{\pm\sqrt{-\Delta}\bigr\}\times Q^4\times \p^2$ is 
birational to an
irreducible component of 
$\sym^{2}\bigl( \OG(\p^2, Q^4)\bigr)$ 
(or a union of 2  irreducible components if $-\Delta$ is a square in $k$)
by (\ref{sq.og.some.say}.4).
Thus Theorem \ref{n>4.main.thm} can be viewed as  a  generalization of  Theorem \ref{4-folds.3.thm}. Note, however that    Theorem \ref{n>4.main.thm}
fails if  $Q^n(k)= \emptyset$. Indeed, for  $n\geq 5$ consider the diagonal quadric
$Q^n:=(x_0^2+\cdots+ x_{n+1}^2=0)$ over any field $k\subset \r$.
Then $  \maps_2(\p^2, Q^n)(k)= \emptyset$ but the conjugate pair of
linear spaces  
$$
L^2_{\pm}:=(x_0\pm ix_1=x_2\pm ix_3=x_4\pm ix_5=x_6=\cdots=x_{n+1}=0)
$$
shows that $ \sym^2\bigl(\OG(\p^2, Q^n)\neq  \emptyset$.

The birational equivalence in Theorem \ref{n>4.main.thm} is not
completely $\aut(Q^n)$-equivariant. The reason is that we usually imagine the left hand side as $\chow_V(Q^n)\times \aut(\p^2)$. However, the proof works with 
pointed Veronese surfaces. Thus it views the left had side as $\univ_V(Q^n)\times \aut(p\in \p^2)$, where $\univ_V(Q^n)\to \chow_V(Q^n)$
denotes the universal family of Veronese surfaces in $Q^n$
and $\aut(p\in \p^2)$ is the group of automorphisms of $\p^2$ that fix the marked point $p$.
We have $\aut(Q^n)$-equivariance for the latter presentation.

Theorem \ref{n>4.main.thm} implies that
$\maps_2(\p^2, Q^n) $  is geometrically irreducible for $n\geq 5$. In particular, 
  projected Veronese surfaces  and quadruple planes can be deformed to Veronese surfaces in a 5-dimensional quadric;
see Example \ref{n>4.main.thm.exmp.2} for explicit descriptions.
\end{rem}

I do not know a  birational classification of 
 symmetric squares of orthogonal Grassmannians in general.
We give an  explicit description of  
$\sym^2\bigl(\OG(\p^1, Q^3)\bigr)$ in (\ref{sym2.OG.Q3.prop.pf}) and this implies the following.

\begin{prop}\label{sym2.OG.Q3.prop}
 Let $Q^5=\bigl(q(x_0,\dots, x_4)=x_5x_6\bigr)$ be a smooth quadric 5-fold
over  a field $k$. Then
$$
\maps_2(\p^2, Q^5)
\simb 
\bigl(\Delta(q)z^2=q(x_0,\dots, x_4)\bigr)\times \p^{22},
$$
where $\Delta(q)$ is the discriminant of $q$.
\end{prop}

Note that  the quadratic form 
$\Delta(q)z^2-q(x_0,\dots, x_4)$ is never positive definite, thus the
Hasse-Minkowski theorem (cf.\ \cite[p.41]{serre-ca}) and Theorem \ref{n>4.main.thm} imply the following.

\begin{cor}  Let  $Q^5$ be  a smooth quadric 5-fold over $k$
such that   $Q^5(k)\neq \emptyset$.
Assume that $k$ is one of the following:
\begin{enumerate}
\item $\r$ or any real closed field,
\item a number field,  
\item  a function field of a curve over a finite field or 
\item a function field of a surface
over an algebraically closed field. 
\end{enumerate}
Then $\maps_2(\p^2, Q^5)\simb \p^{26}$. \qed
\end{cor}

For other fields the situation is more complicated.

\begin{exmp} Let $K:=k(t_1,\dots, t_4)$ where the $t_i$ are 
algebraically independent. Consider the   quadric 5-fold
$Q^5:=(x_0^2+t_1x_1^2+\cdots+t_4x_4^2=x_5x_6)$.
It is easy to see that
$$
x_0^2+t_1x_1^2+\cdots+t_4x_4^2=t_1t_2t_3t_4\cdot z^2
$$
has no solutions in $K$. Therefore, using Theorem \ref{n>4.main.thm} and
Proposition \ref{sym2.OG.Q3.prop}, we obtain that 
 $\sym^2\bigl(\OG(\p^1, Q^{3}_W)\bigr)(K)=\emptyset $ and hence
$\maps_2^{\circ}(\p^2, Q^5)(K)=\emptyset  $. 

Thus the image of any  quadratic map  $\p^2\map Q^5$ is a point or a conic.
Examples for the latter are
$(u{:}v{:}w)\mapsto (0:\cdots :0: uv: u^2: t_4v^2)$. These correspond to singular $K$-points of $\maps_2(\p^2, Q^n) $. 
\end{exmp}

The question of uniqueness  of  Veronese surfaces
$V\subset Q^n$
(up to coordinate changes by $\aut(\p^2)\times\aut(Q^n)$)
is not as simple as in dimension 4.

\begin{prop} \label{V.UP.TO.ISOM.5.THM}
Let $k$ be a field. A smooth quadric hypersurface $Q^5\subset \p^6$ contains a Veronese surface iff its equation can be written as
$$
Q^5_{a,b,c}=\bigl(x_0x_5=x_1^2+ax_2^2+bx_3^2+abx_4^4+cx_6^2\bigr).
$$
Up to coordinate changes by $\aut(\p^2)\times\aut\bigl(Q^5_{a,b,c}\bigr)$,
the quadric $Q^5_{a,b,c}$ contains a unique Veronese surface. 
A  representative  is given by
$$
(u{:}v{:}w)\mapsto 
\bigl(u^2+aw^2:
uv : vw : uw : w^2 : v^2+bw^2: 0\bigr).
$$
\end{prop}

The proof is given in (\ref{pf.V.UP.TO.ISOM.5.THM}) and the next examples show
  that we have  non-uniqueness for $n\geq 6$; see also 
Proposition \ref{real.sign.prop}.

\begin{exmp} \label{nonun.exmp}
If $k$ is algebraically closed then all smooth subquadrics
$Q^4\subset Q^n$ are isomorphic, thus Theorem \ref{V.UP.TO.ISOM.THM}
implies that there is a unique Veronese surface
$V\subset Q^n$,
up to coordinate changes by $\aut(\p^2)\times\aut(Q^n)$.
The same holds over $\r$ if $Q^n$ is the sphere.

However, for most other fields there are  quadrics that contain
inequivalent Veronese surfaces,  up to coordinate changes.
Indeed, by Theorem \ref{V.UP.TO.ISOM.THM}, isomorphism classes of quadric 4-folds that contain a Veronese surface are in one-to-one correspondence with 
 isomorphism classes of plane conics. Assume that we have
two non-isomorphic ternary forms
$x^2+ay^2+bz^2$ and $x^2+ay^2+cz^2$. Consider the quadric 6-fold
$$
Q^6_{a,b,c}=\bigl(x_0x_7=x_1^2+ax_2^2+bx_3^2+abx_4^4+cx_5^2+acx_6^4\bigr).
\eqno{(\ref{nonun.exmp}.1)}
$$
It contains a  Veronese surface $V_b$  that satisfies
$(x_5= x_6=0)$ and also a $V_c$  that satisfies
$(x_3= x_4=0)$. These are inequivalent since the quadrics they span are not isomorphic.

As another example,  let $k$ be a number field  and consider
$$
Q^7:=\bigl(x_0^2+x_1x_2+x_3x_4+x_5x_6+x_7x_8=0\bigr).
\eqno{(\ref{nonun.exmp}.2)}
$$
Note that $au^2-av^2\sim uv$, so every $Q^4_{a,b}$ has at least one  embedding into $Q^7$. Thus there are infinitely many inequivalent 
 embeddings  $V\into Q^7$ 
 defined over $k$,  up to coordinate changes.

\end{exmp}

\begin{notation} \label{basic.not.not}
From now on, $Q^n:=\bigl(q(x_0,\dots, x_{n+1})=0 \bigr)\subset \p^{n+1}$ denotes a smooth quadric hypersurface over a field $k$. 
Starting with Section \ref{sec.degmaps} we assume that $\chr k\neq 2$. 
The explicit formulas are usually worked out only for diagonal forms.

More generally, $Q^{n,i}$  denotes a  quadric hypersurface over a field $k$
whose non-smooth locus has dimension $i$. We use this mainly for $i=0$. 

A quadratic form $q$ over $k$---or the corresponding quadric  hypersurface $Q$ over $k$---is called
{\it isotropic} if $Q(k)\neq \emptyset$ and {\it anisotropic} if $Q(k)= \emptyset$.

The {\it discriminant} of a symmetric bilinear  form 
$q=\sum_{0\leq i, j\leq n+1}a_{ij}x_ix_j$ over $k$ is 
$\Delta(q):=\det (a_{ij})\in k/(k^*)^{2}$. 
This defines the  discriminant of any quadratic form whenever $\chr k\neq 2$. (Note that some authors use different signs and 2-powers.)
A quadric $Q^n=\bigl(q(x_0,\dots, x_{n+1})=0 \bigr)$ determines the corresponding quadratic form $q$ up to a multiplicative scalar $\lambda$.
This changes the discriminant by $\lambda ^{n+2}$. So
the {\it discriminant} of $Q^n$ is only defined  for even dimensional quadrics and then it is
$$
\Delta(Q^n):=\Delta\bigl(q(x_0,\dots, x_{n+1})\bigr)\in k/(k^*)^{2}.
$$
Note that $Q^n$ is smooth iff $\Delta(Q^n)\neq 0$.

For any point  $p\in Q^n$, let $Q^{n-1,0}_p\subset Q^n$ denote the intersection of $Q^n$ with its tangent plane at $p$. Thus
$Q^{n-1,0}_p$ is an $(n-1)$-dimensional quadric cone 
over an $(n-2)$-dimensional smooth quadric---the projectivized tangent cone---denoted by $Q^{n-2}_p$.
Both $Q^{n-1,0}_p$ and $Q^{n-2}_p$ are 
defined over the field $k(p)$.  The  $Q^{n-2}_p$ together form the universal projectivized tangent cone $\PTC(Q^n)\to Q^n$.

Assume next that  $Q^n(k)\neq \emptyset$. Pick any $p\in Q^n(k)$.
In suitable coordinates we may assume that $p=(0,\dots, 0,1)$ and 
$$
q(x_0,\dots, x_{n+1})\sim q_p(x_0,\dots, x_{n-1})+x_nx_{n+1}.
$$
Then $Q^{n-2}_p\cong \bigl(q_p(x_0,\dots, x_{n-1})=0\bigr)\subset \p^{n-1}$.
By Theorem \ref{witt.cancel}, $q_p$ and  $Q^{n-2}_p$ are independent of $p$;
we denote this quadric by $Q^{n-2}_W$ and call it the
{\it Witt reduction} of $Q^n$.  We can iterate the  Witt reduction
until we reach an anisotropic  quadric
$$
Q_{\rm an} =Q^{n-2r}_{\rm an} = \bigl(q_{\rm an}(x_0,\dots, x_{n+1-2r})=0\bigr),
$$
 called the {\it anisotropic kernel} of $Q^n$ (or of $q$). The number 
$r$ is called  the {\it Witt index.} 
Equivalently, $r-1$ is the maximum dimension of a linear space that is
defined over $k$ and
contained in $Q^n$. 

A form is called {\it split} if  $2r=n+2$,
equivalently, if  it can be written as
$$
x_0x_1+x_2x_3+\cdots+x_{2r-2}x_{2r-1}.
$$
Let $\pi_p:Q^n\map \p^n$ denote the projection from $p$.
We can factor it as the blow-up of $p$ followed by the contraction of the
birational transform $E^{n-1,0}_p$ of $Q^{n-1,0}_p$. Here $E^{n-1,0}_p$ is a 
$\p^1$-bundle over $Q^{n-2}_p\subset \p^n$. This shows that there is an equivalence
$$
\bigl\{p\in  Q^n\bigr\}
\ \leftrightarrow\
\bigl\{Q^{n-2}_W\subset \p^n\bigr\}
\eqno{(\ref{basic.not.not}.1)}
$$
between $n$-dimensional pointed quadrics and 
$(n-2)$-dimensional  quadrics in $\p^n$.

A geometric proof of the special case of Witt's theorem that we used  above is the following. 
 Fix a point  $p_0\in Q^n(k)$. For any $p\in Q^n\setminus Q^{n-1,0}_{p_0}$ we get natural isomorphisms  (defined over $k(p)$)
$$
Q^{n-2}_p\cong Q^{n-1,0}_p\cap Q^{n-1,0}_{p_0}\cong Q^{n-2}_{p_0}.
\eqno{(\ref{basic.not.not}.2)}
$$
This also shows that  the universal projectivized tangent cone $\PTC(Q^n)\to Q^n$ is birationally trivial
$$
\PTC(Q^n)\simb Q^n\times Q^{n-2}_W.
\eqno{(\ref{basic.not.not}.3)}
$$
\end{notation}

\begin{say}[Witt's cancellation theorem]\label{witt.cancel}
Let $V$ be a vector space and $q$ a non-degenerate quadratic form on $V$. 
Let $V_1, V_2\subset V$ be vector subspaces such that $\bigl(V_i, q|_{V_i}\bigr)$
are non-degenerate and 
there is an isometry  $\phi:\bigl(V_1, q|_{V_1}\bigr)\to \bigl(V_2, q|_{V_2}\bigr)$.
One form of Witt's cancellation theorem says that  $\phi$ extends to an isometry  $\Phi: V\to V$.

This is one example where the correspondence between quadrics and quadratic forms gets cumbersome.
For example, consider  $(\r^5, q:=x_1^2+ x_2^2+ x_3^2- x_4^2- x_5^2)$,
 $V_1:=(x_3=x_5=0)$ and $V_2:=(x_1=x_3=0)$. Note that
the quadrics  $(q=x_3=x_5=0)$ and $(q=x_1=x_3=0)$ are isomorphic but
the quadratic forms $\bigl(V_1, q|_{V_1}\bigr)$ and $ \bigl(V_2, q|_{V_2}\bigr)$ are not. There is no automorphism of the quadric $(q=0)$ that carries
 $(q=x_3=x_5=0)$ to  $(q=x_1=x_3=0)$. 
\end{say}

\begin{say}[Effective parametrization] The methods of this paper are effective
for 4-dimensional quadrics, that is, the proofs can be converted into explicit formulas
giving all elements of $\maps_2(\p^2, Q^4)(k)$ for any field $k$, provided we know the answer to some basic questions about $Q^4$.

If there is a $k$-map  $\phi:\p^2\map Q^n$ then
$Q^n(k)\neq \emptyset$. Conversely, 
 if $p\in Q^n(k)$ is any $k$-point then the constant map $\p^2\to \{p\}\into Q^n$
is in $\maps_2(\p^2, Q^n)(k)$. Thus $\maps_2(\p^2, Q^n)(k)\neq \emptyset$
iff $Q^n(k)\neq \emptyset$. To the best of my knowledge, there is no algorithm that decide whether a quadric over an arbitrary field has a $k$-point or not,
and we do not claim to say anything about this question. 

We show, however, that once we can decide the existence of points, we can also effectively write down all maps in $\maps_2(\p^2, Q^n)(k)$.
To be precise, we need to decide the  existence of points also for certain subforms and the  condition we need is the following.
We start with the cases when $n\leq 5$.
\medskip

{\it Assumption.} For any quadratic form $q(x_0,\dots, x_{n+1})$
 we are given (or are able to find) a Witt decomposition of $q$. That is, a  linear change of coordinates
$$
q(x_0,\dots, x_{n+1})\sim q_{\rm an}(y_0,\dots, y_{n-2r+1})+y_{n-2r+2}y_{n-2r+3}+\cdots+
y_{n}y_{n+1}
$$
where $q_{\rm an}$ is anisotropic.

Note that such a linear change of coordinates always exists over any field, but
in general we do not know how to find it. Given a field $k$,
there is an algorithm to find such a coordinate change for every quadratic form over $k$ iff there is an algorithm that for every  quadratic form $q(x_0,\dots, x_{r+1})$ over $k$

-- either proves that $q$  is anisotropic,

-- or finds a non-trivial solution  $q(a_0,\dots, a_{r+1})=0$. 

Thus the assumption is always satisfied of $k$ is algebraically closed or
real closed. Over  any number field, the Hasse--Minkowski theorem
(see, for instance  \cite[p.41]{serre-ca}) provides a (very inefficient)  algorithm.
\medskip

For $n\geq 6$ we also need to know a parametrization of degree 2 points of
the orthogonal Grassmannian $\OG(\p^2, Q^n)$. I do not even know how to do this over nice fields like $\q$ or $\c(t)$.

\end{say}

\section{Surfaces containing many circles}

I got interested in these questions after reading
the  papers \cite{2015arXiv150306481S, skopenkov} which 
furnish the last step of a project, started by Kummer and Darboux,  
to describe all surfaces in $\r^n$ that contain at least 2 circles through  every point;
see also \cite{MR1824751, MR2871363}.

The inverse of stereographic projection connects this problem with 
real algebraic surfaces $F$ on the sphere 
$$
 \s^n:= (x_1^2+\dots+x_{n+1}^2=x_0^2)\subset \r\p^{n+1}.
$$
The formulas for stereographic projection  are nicest if we project
the sphere
$$
\s^n:=\bigl(x_1^2+\cdots+x_n^2+z^2=1\bigr)
$$
from the south pole  $(0,\dots, 0,-1)$.
Then
$$
\pi(x_1,\dots,x_n, z)=\left(
\frac{x_1}{1+z},\dots, \frac{x_n}{1+z}\right)
$$
with inverse
$$
\pi^{-1}(x_1,\dots,x_n)=
\left(
\frac{2x_1}{1+\Sigma},\dots, \frac{2x_n}{1+\Sigma},\frac{1-\Sigma}{1+\Sigma} \right),
$$
where $\Sigma=\sum_i x_i^2$. 
A theorem that  Ptolemy attributes to Hipparchus  ($\sim$ 190-120 BC)
says that stereographic projection preserves circles;   a modern treatment was given by Halley in 1695.

Any conic on a sphere is a circle, thus 
describing all surfaces in $\r^n$ that contain many circles
is equivalent to 
 describing real algebraic surfaces
$$
F\subset \s^n:= (x_1^2+\dots+x_{n+1}^2=x_0^2)\subset \r\p^{n+1}
$$
that contain many conics.

Building on works of \cite{kachi-99} and \cite{MR1824751}, 
the hardest cases are when $F$ is a projection of
the Veronese surface or of a degree 8 del~Pezzo surface.
These are classified in \cite{skopenkov}
for  $\s^4$.
His idea is to  rewrite the equation as
$$
 x_1^2+x_2^2+x_3^2+x_4^2=(x_0-x_5)(x_0+x_5),
$$
  think of
$x_1^2+x_2^2+x_3^2+x_4^2$ as the norm of the quaternion
${\mathbf X}:=x_1+x_2i+x_3j+x_4k$ and solve the quaternionic equation
$$
\norm ({\mathbf X})=z_0z_5,
\qtq{where} {\mathbf X}\in {\mathbb{H}}[u,v],\ z_0, z_5\in {\mathbb{R}}[u,v].
$$
 This is quite nontrivial since
${\mathbb{H}}[u,v]$ is far from having  unique factorization.

In our language, we work with the quadratic form
$q:=x_1^2+\dots+x_5^2-x_0^2$ over $\r$. 
Its  discriminant is  $\Delta=-1$ 
and the anisotropic kernel is  $x_1^2+x_2^2+x_3^2+x_4^2$. 
Thus we see from Propositions  \ref{4-folds.2.thm}--\ref{4-folds.1.thm}
and  Theorem \ref{4-folds.3.thm} that the Veronese
components  
$\maps_{V}(\p^2, Q^4)(\r)\simb \r\p^{20}\amalg \r\p^{20}$
give all real solutions. By Theorem \ref{V.UP.TO.ISOM.THM}, there is only one
solution up to isomorphism. 
Since a Veronese surface in $\p^{n+1}$ is contained in a linear subspace
of dimension $\leq 5$, these imply the following.

\begin{thm}\label{sur.circ.all.thm} There are precisely 2 surfaces  $F\subset \s^n$  
(up to  the group of M\"obius transformations $\aut(\s^n)=\OO(x_1^2+\dots+x_{n+1}^2-x_0^2)$) that contain a circle through any 2  points  $p_1, p_2\in F$:
\begin{enumerate}
\item the sphere  $\s^2\subset \s^n$ and
\item the Veronese surface  $\r\p^2\into \s^4\subset \s^n$. \qed
\end{enumerate}
\end{thm}

There are several ways to describe the  real Veronese surface  $\r\p^2\into \s^4$.
For example, one can think of the Veronese surface $V$ as the set of rank 1 matrices in  the projectivization of the space of $3\times 3$ symmetric matrices:
$$
(u:v:w)\mapsto
\left(
\begin{array}{ccc}
u^2 & uv & uw \\
uv & v^2 & vw \\
uw & vw & w^2
\end{array}
\right).
\eqno{(\ref{sur.circ.all.thm}.3)}
$$
The identity
$(uv)^2+(vw)^2+(uw)^2=u^2v^2+v^2w^2+u^2w^2$
becomes the quadratic equation
$$
x_{01}^2+x_{12}^2+x_{02}^2=x_{00}x_{11}+x_{11}x_{22}+x_{00}x_{22},
\eqno{(\ref{sur.circ.all.thm}.4)}
$$
which defines an ellipsoid that contains $V$. 
The right hand side can be diagonalized in a less symmetric form as
$$
\tfrac13\bigl(x_{00}+x_{11}+x_{22}\bigr)^2
-\tfrac1{4}\bigl(x_{00}-x_{11}\bigr)^2
-\tfrac1{12}\bigl(x_{00}+x_{11}-2x_{22}\bigr)^2.
\eqno{(\ref{sur.circ.all.thm}.5)}
$$
This gives an embedding of $\r\p^2$ into the standard sphere of radius 1
$$
(u:v:w)\mapsto \tfrac{\sqrt{3}}{u^2+v^2+w^2}
\Bigl( uv: vw : uw: \tfrac{1}{2}(u^2-v^2) : \tfrac{1}{\sqrt{12}} (u^2+v^2-2w^2)\Bigr).
\eqno{(\ref{sur.circ.all.thm}.6)}
$$
Another form is given in Remark \ref{quat.form.rem}:
$$
(u:v:w)\mapsto \tfrac{2}{u^2+v^2+2w^2}
\Bigl( uv: vw : uw: w^2 : \tfrac{1}{2} (u^2-v^2)\Bigr).
\eqno{(\ref{sur.circ.all.thm}.7)}
$$
Another form of the  Veronese embedding into $\r^5$ is obtained by  looking at the  vector space $ V^{(3,3)}$ of $3\times 3$ symmetric matrices of trace 0. $\SO(3,\r)$ acts on $V^{(3,3)}$ by conjugation and the action preserves the eigenvalues. Furthermore, the orbits are in one-to-one correspondence with the ordered set of eigenvalues  $\lambda_1\leq \lambda_2\leq \lambda_3$. 
The generic orbit is large  (with stabilizer  $(\z/2)^2$)
but if $2$ eigenvalues coincide, say for
$-1, -1, 2$, then the orbit is  isomorphic to the Veronese embedding of $\r\p^2$.  The $\SO(3,\r)$-invariant sphere containing this Veronese surface is given by
the equation
$\tr(M^2)=6$.

Note that  $M\mapsto M+{\bf 1}_3$ sends  a
symmetric matrix with eigenvalues $-1, -1, 2$ to a 
symmetric matrix of rank 1  and  trace $3$. This shows that the
above affine form is equivalent to the projective one.
\medskip

Other real quadric hypersurfaces  can contain more then one
families of Veronese surfaces.

\begin{prop} \label{real.sign.prop}
Let $Q^n\subset \r\p^{n+1}$ be a quadric of dimension $n\geq 4$ and of signature  $(r,s)$ where $r\geq s$.  Then 
the number of  different Veronese embeddings  $\r\p^2\into Q^n$ 
 (up to coordinate changes by $\aut(\r\p^2)\times \aut(Q^n)$) is
\begin{enumerate}
\item zero if $s=0$ or $r\leq 4, s\leq 2$,
\item one if $r\geq 5$ and $s\in\{1,2\}$ or if $(r,s)=(3,3), (4,3), (4,4)$ and 
\item at most 3 in general. 
\end{enumerate}
\end{prop}

Proof. If $n=4$ then $Q^n$ contains a Veronese surface only when the signature is $(5,1), (1,5)$ or $(3,3)$. Thus we need to enumerate the different solutions to  $(r,s)=(r',s')+(r'',s'')$ where  $(r',s')$ is one of
 $(5,1), (1,5)$ or $(3,3)$. There are clearly at most 3 solutions
but fewer is $s$ is small or for occasional symmetries. \qed
\medskip

It would be quite interesting to understand higher dimensional generalizations of our results.  The following form seems to be the most promising.

\begin{prob}  Classify real algebraic varieties $X\subset \s^n$ such that
 every point pair  $x_1, x_2\in X$  is contained in a circle
$ C(x_1, x_2)\subset X$.
\end{prob}

The paper \cite{ion-rus} classifies  
complex algebraic varieties $X\subset \p^{n+1}$ that are 
{\it conic-connected,} that is, 
 every point pair  $x_1, x_2\in X$  is contained in at least one  conic
$ C(x_1, x_2)\subset X$; see also \cite[Chap.4]{rus-book}.
Thus one  needs to decide  which conic-connected varieties give real examples in $\s^n$.
There are many  conic-connected  Fano varieties of Picard number 1 and these are not yet fully understood. However,  if $\dim X\leq 6$, then these are all so called  {\it  Mukai manifolds} and a complete
 classification of Mukai manifolds is given in \cite{MR995400, MR1675146}.
Thus   the list of \cite{ion-rus}  is complete for $\dim X\leq 6$, 
 hence it could be  feasible to answer
the above problem when $\dim X\leq 6$.

An interesting extremal case is the  Veronese embedding of $\p^n$ by quadrics.
By a theorem of Gallarati (see \cite[3.4.4]{rus-book})
if $X^n\subset \p^N$ is a linearly non-degenerate   conic-connected variety
and  $N\geq \frac12n(n+3)$ then   $N= \frac12n(n+3)$
and $X^n$ is the  Veronese embedding of $\p^n$ by quadrics.
See \cite{MR0293645} for a version for differentiable manifolds.
One should study which quadrics contain these higher dimensional Veronese varieties.

\section{Maps of $\p^1$ to quadrics}\label{sec.2}

Here we describe the spaces $\maps_d(\p^1, Q^n)$
parametrizing degree $d$ maps from $\p^1$ to a quadric. 
Equivalently, we describe  
 solutions of a homogeneous
quadratic equation $q(x_0,\dots, x_{n+1})=0$, defined over a field $k$,  where the
$x_i=h_i(u,v)$ are  degree $d$   homogeneous polynomials in  2 variables $u,v$.

The proof of Theorem \ref{curve.map.thm} is a direct consequence of
the following 2 propositions.

\begin{prop} \label{maps.to.Q.thm} Let $Q^n$ be a smooth quadric of dimension $n\geq 3$. Then 
$$
\maps_{d}(\p^1, Q^n)\simb \maps_{d-2}(\p^1, Q^n)\times \p^{2n}
\qtq{for $d\geq 3$.}
$$
\end{prop}

\begin{prop} \label{maps.to.Q.d=12} Let $Q^n$ be a smooth quadric of dimension $n\geq 3$. Then 
$$
\maps_{1}(\p^1, Q^n)\simb \OG(\p^1, Q^n)\times \p^3 \qtq{and}
\maps_2(\p^1, Q^n)\simb Q^n\times \p^{2n}.
$$
\end{prop}

Proof. The image of a degree 1 map is a line, this gives
the formula for $ \maps_{1}(\p^1, Q^n)$.

Let $\maps_{2}^{\circ}(\p^1, Q^n) $ denote those degree 2 maps that
map $\p^1$ isomorphically onto a conic.
By \cite[II.3.14]{rc-book},  $\maps_{2}^{\circ}(\p^1, Q^n) $ is an open and dense  subset of $\maps_{2}(\p^1, Q^n) $; this can also be seen by a simple dimension count here.
A conic  is uniquely determined by any 3 of its points, say the images of $0,1,\infty\in \p^1$. This gives an open embedding (hence a birational equivalence)
$$
\maps_{2}^{\circ}(\p^1, Q^n)\into Q^n\times Q^n\times Q^n.
$$
The birational equivalence $ Q^n\times Q^n\times Q^n \simb Q^n\times \p^{2n}$
can be seen as follows. 
Fix two hyperplanes $H_0, H_1\subset \p^{n+1}$. 
Given $p\in Q^n$ and $q_i\in H_i$, there is a unique
degree 2 map $\phi:\p^1\to Q^n$ such that $\phi(\infty)=p$ and, 
$\phi(i)$ is the residual intersection point of $\langle p, q_i\rangle\cap Q^n$ for $i=0,1$.\qed
\medskip

Note that $n=2$ is exceptional for Proposition \ref{maps.to.Q.d=12} and
$$
\maps_2(\p^1, Q^2)\simb \bigl(Q^2\times \p^{4}\bigr)\amalg 
\bigl(\OG(\p^1, Q^2)\times \p^5\bigr).
$$

The plan to prove Proposition \ref{maps.to.Q.thm} is the following.
Given a degree $d$ map $\phi:\p^1\to Q^n$, we aim to write its image as the directrix of a ruled surface $S$. For suitable choice of $S$ the residual intersection of $S$ and $Q^n$ is the image of  a degree $d-2$ map $\psi:\p^1\to Q^n$. We start with a general discussion on ruled surfaces.

\begin{say}[Ruled surfaces]\label{ruled.surf.say}
  Let $C$ be a smooth projective curve and $\phi_i:C\to \p^n$  two morphisms.
We want to understand the ruled surface swept out by the lines
$\langle \phi_1(p), \phi_2(p)\rangle$ for $p\in C$.

Consider the abstract ruled surface  $S:=\p_C\bigl(\phi_1^*\o_{\p^n}(1)+\phi_2^*\o_{\p^n}(1)\bigr)$. 
It has 2 natural disjoint sections  $C_1, C_2\subset S$ 
and $\o_S(1)|_{C_i}\cong \phi_i^*\o_{\p^n}(1)$. Consider the exact sequence
$$
0\to \o_S(1)(-C_1-C_2)\to \o_S(1)\to \o_S(1)|_{C_1}+\o_S(1)|_{C_2}\to 0.
$$
Since $\o_S(1)(-C_1-C_2)$ has degree $-1$ on the rulings, all of its cohomologies are 0. Thus we get that
$$
H^0\bigl(S, \o_S(1)\bigr)\cong H^0\bigl(C, \phi_1^*\o_{\p^n}(1)\bigr)+
H^0\bigl(C, \phi_2^*\o_{\p^n}(1)\bigr).
$$
That is, the pair of morphisms $(\phi_1, \phi_2)$ 
uniquely extends to a morphism
$\Phi:S\to \p^n$. The degree of the image, that is,
the self-intersection of $\Phi^*\o_{\p^n}(1)$ is $\deg \phi_1+\deg\phi_2$.

We can do better if there are $r$ points  $p_1,\dots, p_r\in C$
such that  $\phi_1(p_i)=\phi_2(p_i)$ holds. 
More generally, let $Z\subset C$ be a 0-dimensional subscheme
such that $\phi_1|_Z=\phi_2|_Z$. Let us denote this map by $\phi_Z$. 
We have natural restriction maps
$r_i: \phi_i^*\o_{\p^n}(1)\to \phi_Z^*\o_{\p^n}(1)$. 
We define the sheaf $E$ as the kernel of the map
$$
r_1-r_2: \phi_1^*\o_{\p^n}(1)\oplus \phi_2^*\o_{\p^n}(1)
\to \phi_Z^*\o_{\p^n}(1).
$$
Then $E$ is locally free of rank 2 and we have natural surjections
$E\onto  \phi_i^*\o_{\p^n}(1)$.  Now we set $S:=\p_C(E)$.
As before it has 2 natural  sections  $C_1, C_2\subset S$  such that
 $\o_S(1)|_{C_i}\cong \phi_i^*\o_{\p^n}(1)$, but now
$C_1\cap C_2\cong Z$. The pair of morphisms $(\phi_1, \phi_2)$ 
uniquely extends to a morphism $\phi_{12}: C_1\cup C_2\to \p^n$.
As before, there is a unique extension to 
$\Phi:S\to \p^n$. The gain is that   
the self-intersection of $\Phi^*\o_{\p^n}(1)$ is $\deg \phi_1+\deg\phi_2-\deg Z$.

Let now $Q^n\subset \p^{n+1}$ be a smooth quadric, $C$ a smooth curve  and $\phi:C\to Q^n$
a morphism of degree $a$. Choose any other morphism $\psi:C\to \p^{n+1}$
of degree $b$. We get a ruled surface
$\Phi:S(\phi,\psi)\to \p^{n+1}$. Assume that its image is not contained in $Q^n$.
Then $\Phi^{-1}(Q^n)$ is in the linear system of sections of $\o_S(2)=\Phi^*\o_{\p^{n+1}}(2)$. By construction it contains
$C_1$; let $R:=R(\phi,\psi)$ denote  the residual intersection.
Thus $R$ is a union of a section $C_3\cong C$ of $S$ plus possibly some rulings.
These rulings are  over those points $p\in C$ for which
the line $\langle \phi(p), \psi(p)\rangle$ is contained in $Q^n$
and possibly over those points $p\in C$ for which
 $ \phi(p)= \psi(p)$.  Thus we get a new morphism
$$
\phi*\psi:\Phi|_{C_3}:C\to Q^n.
\eqno{(\ref{ruled.surf.say}.1)}
$$  
If there is  a 0-dimensional subscheme $Z\subset C$ 
such that $\phi|_Z=\psi|_Z$ then 
the degree of the residual curve $R$ is  $2\deg S-\deg \phi=\deg \phi+2\deg \psi-2\deg Z$ hence
$$
\deg (  \phi*\psi)= \deg \phi+2\deg \psi-2\deg Z-\#(\mbox{rulings in $R$}).
\eqno{(\ref{ruled.surf.say}.2)}
$$
We get an interesting degenerate case if $\psi$ also maps to $Q^n$.
Then $\Phi^{-1}(Q^n)$ contains both sections $C_1, C_2$, hence
the residual curve $R:=R(\phi,\psi)$ is a union of 
$\deg \phi+\deg \psi$ rulings.
\end{say}

\begin{say}[Proof of (\ref{maps.to.Q.thm})] \label{pf.of.maps.to.Q.thm}
We aim to prove the  birational equivalence
$$
\maps_{d}(\p^1, Q^n)\simb \maps_{d-2}(\p^1, Q^n)\times \p^{2n}
\qtq{for $d\geq 3$.}
\eqno{(\ref{pf.of.maps.to.Q.thm}.1)}
$$
We fix the points   $0,1,\infty\in \p^1$, these will play a special role
in the construction.  To emphasize this, we write
$\maps_d(\p^1, 0,1,\infty; Q^n) $ instead of $\maps_{d}(\p^1, Q^n) $,
though these spaces are isomorphic.

 Fix an auxiliary hyperplane $L^n\subset \p^{n+1}$.  We construct a morphism
$$
\Pi_d: \maps_d(\p^1, 0,1,\infty; Q^n)\map \maps_{d-2}(\p^1, 0,1,\infty; Q^n).
$$
Let $\phi_d: (\p^1, 0,1,\infty)\to Q^n$ be a morphism of degree $d$.
Let $\phi_1: (\p^1, 0,1,\infty)\to \p^{n+1}$ be a morphism of degree $1$
that sends $0\mapsto \phi_d(0), \infty\mapsto \phi_d(\infty)$ and
$1$ to the intersection of $L^n$ with the line
$\langle\phi_d(0),\phi_d(\infty)\rangle$.  As in (\ref{ruled.surf.say})
the pair $(\phi_d, \phi_1)$ determines a ruled surface $S_{d-1}$ and 
by (\ref{ruled.surf.say}.2) $Q^n\cap S_{d-1}=\phi(\p^1)\cup R_{d-2}$
where $R_{d-2}$ is a curve of degree $d-2$. We check below that it is an 
irreducible rational curve for general $\phi_d$.
Since both $\phi_d$ and  $\psi_{d-2}$ give sections of $S_{d-1}$,
the 3 marked points on the image of  $\phi_d$ give
3 marked points on the image of $\psi_{d-2}$. 
This determines  $\psi_{d-2}:(\p^1, 0,1,\infty)\to Q^n$.

Next we show that the generic fiber of $\Pi_d$ is birational to
$\p^{2n}$. To see this, fix $\psi_{d-2}:(\p^1, 0,1,\infty)\to Q^n$.
Pick points  $x_0, x_{\infty}\in Q^n$. Let 
 $\phi_1: (\p^1, 0,1,\infty)\to \p^{n+1}$ be a morphism of degree $1$
that sends $0\mapsto x_0, \infty\mapsto x_{\infty}$ and
$1$ to the intersection of $L^n$ with the line connecting 
$x_0, x_{\infty}$.  As in (\ref{ruled.surf.say})
the pair $(\psi_{d-2}, \psi_1)$ determines a ruled surface $S_{d-1}$.
By construction $Q^n\cap S_{d-1}=\psi(\p^1)\cup R_{d}$. 
As before, this determines  $\phi_{d}:(\p^1, 0,1,\infty)\to Q^n$. 
We need to check that the curves $R_{d-2}$ and $R_d$ are irreducible for
general choices of $\phi$ and $\psi$.
For this it is enough to find one particular case where the second
construction
 gives an irreducible curve $R_d$; the same set-up then proves the converse too. We start with the degenerate case
where $\psi_{d-2}$ maps onto a line  $L\subset Q^n$.

We can now vary $Q^n$ in the linear system  $|\o_{\p^{n+1}}(2)|(-L)$,
which is base-point free outside $L$. 
Bertini theorem gives that for general $Q^n$ 
the residual intersection  $R_d:=(S_{d-1}\cap Q^n)\setminus L$ 
is irreducible. Since all smooth quadrics of dimension $n$ are
isomorphic over $\bar k$, irreducibility can be achieved by fixing $Q^n$ and 
changing $S_{d-1}$.
\end{say}

Note that if $n=2$ and $\phi_d$ is a curve of bidegree $(1,d-1)$ on $Q^2$ then
the resulting  $R_{d-1}$ is a union of $d-1$ lines. So one really needs to check
that the construction gives an irreducible curve for $n\geq 3$.

\section{Rational curves on quadrics}\label{sec.3}

Instead of looking at degree $d$  maps from $\p^1$ to $Q^n$, it is also of interest to study  geometrically rational curves contained in $Q^n$.

\begin{defn} \label{mo.defn}
Let $X$ be a projective variety and $L$ an ample line bundle on $X$. For $g, d\geq 0$ let  $\mo_g(X,d)$ be the space parametrizing morphisms
$\phi:C\to X$ where

-- $C$ is a smooth, irreducible curve of genus $g$,

-- $\phi$ is a morphism such that $C\to \phi(C)$ is birational and

-- $\deg_C\phi^*L=d$.

One can think of $\mo_g(X,d)$ as a subset of the Chow variety $\chow(X)$ 
and  it is frequently denoted by $\rc_d(X)$  \cite[Sec.II.2]{rc-book}.

One can also view  $\mo_g(X,d)$ 
as an open subset of $\bar M_g(X,d)$,
the moduli space of stable maps of genus $g$ and of degree $d$ \cite{ful-pan}.
For most purposes the latter is the best compactification but
for the birational properties of $\mo_g(X,d)$ the precise
compactification does not matter much.

Note that 
$\maps_{d}(\p^1, X)$ is a $\PGL_2$-torsor over $\mo_0(X,d)$ 
but this torsor is usually not  Zariski-locally trivial.
However,
if $d=2e+1$ is odd and $(\phi:C\to X)$ is in $\mo_0(X,d)$
then $\psi^*L\otimes \omega_{C_n}^e$ is a degree 1  line bundle on $C$.
Thus the universal family over $\mo_0(X,d)$ 
is Zariski-locally trivial and we conclude that 
$$
\maps_{2e+1}(\p^1, X)\simb \mo_0(X,2e+1)\times \PGL_2\simb \mo_0(X,2e+1)\times \p^3.
\eqno{(\ref{mo.defn}.1)}
$$
However, for $d=2e$ even, there can be (and usually there are)
degree $2e$ rational curves on $X$ without $k$-points. Thus we still get a
rational map $\maps_{2e+1}(\p^1, X)\map \mo_0(X,2e+1)$
but a typical fiber is a non-trivial principal homogeneous space under
$\PGL_2$. It seems much harder to connect the
birational properties of  $\maps_{2e+1}(\p^1, X)$ and $ \mo_0(X,2e+1)$.
\end{defn}

We are mostly interested in the spaces $\mo_0(Q^n,d)$.
If $n\geq 3$ then a dense open subset of
 $\mo_0(Q^n,d)$ parametrizes smooth rational curves in $Q^n$;
 see \cite[II.3.14]{rc-book}. 
For odd values of $d$, Theorem \ref{curve.map.thm} and (\ref{mo.defn}.1)
give the following.

\begin{cor} \label{RC.quad.odd.thm}
Let $Q^n$ be a smooth quadric of dimension $n\geq 3$. Then
 $$
\mo_0(Q^n,2e+1)\times \p^3\simb \OG(\p^1, Q^n)\times \p^{2en+3}. \qed
$$
\end{cor}

It is quite likely that the $\p^3$ factor can be canceled. This is not hard to see for $d=3$  but the above argument  does not establish it. 

We are much less successful  for even degrees and for other genera; we get a
complete description in 4 cases only.

\begin{thm} \label{RC.quad.thm}
Let $Q^n$ be a smooth quadric of dimension $n$. Then
\begin{enumerate}
\item $\mo_0(Q^n,2)\simb \p^{3n-3}$ if $n\geq 1$,
\item $\mo_1(Q^n,4)\simb \p^{4n}$  if $n\geq 2$,
\item $\mo_0(Q^n,4)\times \p^2\simb \sym^2\bigl(\OG(\p^1, Q^n)\bigr)\times \p^{n+5}$  if $n\geq 3$.
\item $\mo_1(Q^n,5)\times \p^6\simb \OG(\p^1, Q^n)\times \p^{5n+3}$  if $n\geq 3$.
\end{enumerate}
\end{thm}

Proof.  The first two claims are clear. A general conic in $Q^n$ is the
intersection of $Q^n$ with a 2-plane, so
$$
\mo_0(Q^n,2)\simb \grass(\p^2,\p^{n+1})\simb \p^{3n-3}.
$$
Similarly, a  general degree 4 elliptic curve in $Q^n$ is the
intersection of $Q^n$ with a 3-plane and a quadric in it, so
$$
\mo_1(Q^n,4)\simb \grass(\p^3,\p^{n+1})\times \p^{8}\simb \p^{4n}.
$$

For (3) and (4) we use  degree 4  del~Pezzo surfaces and  start with  the smallest case $n=3$. 
Pick an auxiliary  fixed point $p_0\in \p^4\setminus Q^3$.

Let $C_4\subset Q^3$ be a degree 4 rational normal curve.
Pick a point pair $\{p_1, p_2\}\in \sym^2(C_4)\cong \p^2$ 
and set $C_2:=Q^3\cap\langle p_0, p_1, p_2\rangle$. Note that $C_4\cup C_2$
is a degree 6 curve of arithmetic genus 1. 
The linear system of quadric sections of  $Q^3$ that contain $C_4\cup C_2$ has dimension $13-11=2$.
We compute in 
(\ref{curve.dp4.say}.1) that the base locus of such a pencil is
$C_4+C_2+\ell_1+\ell_2$. This gives a map
$$
\Pi: \mo_1(Q^3,4)\times \p^2\map \sym^2\bigl(\OG(\p^1, Q^3)\bigr).
$$
What are the fibers? Given a pair of lines  
$\{\ell_1, \ell_2\}\in \sym^2\bigl(\OG(\p^1, Q^3)\bigr)$
pick points $r_1\in \ell_1, r_2\in \ell_2$.  Set 
$C_2:=Q^3\cap\langle p_0, r_1, r_2\rangle$. Note that $B:=C_2+\ell_1+\ell_2$
 is a degree 4 curve of arithmetic genus 0.
Using (\ref{curve.dp4.say}.2) we see that 
  $h^0\bigl(Q^3, \o_{Q^3}(2)(-B)\bigr)=5$ and 
 for any 2 sections of $ \o_{Q^3}(2)(-B) $, the residual intersection is a $C_4$.
Thus the fiber of $\Pi$ is given by the choices of $r_1, r_2$ and
the residual $C_4$. The latter correspond to an open subset of 
  $\grass(\p^1, \p^4)$ which is birational to $\p^6$.
Given a pair of lines $\ell_1, \ell_2\in Q^3\subset \p^4$,
we can project them to a $\p^2$ and then the choice of
the points $r_i\in \ell_i$ corresponds to the dual $\p^2$.
This shows that 
$$
\mo_0(Q^3,4)\times \p^2\simb \sym^2\bigl(\OG(\p^1, Q^3)\bigr)\times \p^8.
$$
Next we consider higher dimensional quadrics. Any $C_4\subset Q^n$ spans a $\p^4$.
This defines a map
$$
\mo_0(Q^n,4)\map \grass(\p^4, \p^{n+1}),
$$
whose fiber over an $L^4\in  \grass(\p^4, \p^{n+1})$
is $\mo_0(L^4\cap Q^n,4)$.  We wrote down above that this is stably birational to the universal family of pairs of lines. 
This family parametrizes pairs $\bigl(\{\ell_1, \ell_2\}, L^4\bigr)$
where $\ell_1, \ell_2$ are lines in $L^4\cap Q^n$. We can also
parametrize this by first choosing 
$\{\ell_1, \ell_2\}\in \sym^2\bigl(\OG(\p^1, Q^n)\bigr)$ and then an $L^4$ that contains their span.  This corresponds to
$\sym^2\bigl(\OG(\p^1, Q^n)\bigr)\times \p^{n-3}$.
Putting these together gives that
$$
\sym^2\bigl(\OG(\p^1, Q^n)\bigr)\times \p^{n-3}\times \p^8
\simb  \mo_0(Q^n,4)\times \p^2. 
$$

In order to prove (4) we  again start with  the smallest case $n=3$. 
The linear system of quadric sections of  $Q^3$ that contain $E_5$ has dimension $13-9=4$. By
(\ref{curve.dp4.say}.2), a general pencil of such quadrics gives a
residual intersection which is a degree 3 rational curve. This gives a map
$$
\Pi:  \mo_1(Q^3,5)\times \p^6\map \mo_0(Q^3, 3).
$$
The fiber is given by pencils in $|\o_{Q^3}(2)|(-C_3)$ which is birational to
$\grass(\p^1, \p^7)\simb \p^{12}$. This proves that
 $$
\Pi:  \mo_1(Q^3,5)\times \p^6\simb \mo_0(Q^3,3)\times \p^{12}\simb
\OG(\p^1,Q^3)\times \p^{18},
$$
where the second birational equivalence follows from (\ref{RC.quad.odd.thm}). 
Extending this to higher dimensions works as before since
 a general $E_5\subset Q^n$ spans a $\p^4$.
This defines a map
$$
\mo_1(Q^n,5)\map \grass(\p^4, \p^{n+1})
$$
whose fiber over an $L^4\in  \grass(\p^4, \p^{n+1})$
is $\mo_1(L^4\cap Q^n,5)$. \qed

\begin{say}[Curves on degree 4 del~Pezzo surfaces]\label{curve.dp4.say}
Let $C\subset S_4$ be  a reduced curve of degree $d$ and arithmetic genus $p_a$
on a  degree 4 del~Pezzo surface  $S=S_4$. Since 
$$
H^0\bigl(S, \o_S(2)\bigr)=13\qtq{and} H^0\bigl(C, \o_S(2)|_C\bigr)\leq 2d+1,
$$
we see that  if $2d+1\leq 12$ then  we get another effective  curve $C'\subset S$ such that  $C+C'\sim \o_{S}(2)$.
The adjunction formula  $2p_a(C)-2=C(C-H)$ gives that
$(C^2)=d+2p_a-2$. From this we get that
the degree of $C'$ is $d'=8-d$ and its   arithmetic genus is
$p_a'=4-d+p_a$. We used some special cases of this.

(\ref{curve.dp4.say}.1)  If $d=6$ and $p_a=1$ then $d'=2$ and $p_a'=-1$. Thus $C'$ is the disjoint union of 2 lines. Conversely, if $d=2$ and $p_a=-1$ then  $|C'|$ is a 6-dimensional linear system whose general member is a smooth elliptic curve.

(\ref{curve.dp4.say}.2)  If $d=5$ and $p_a=1$ then $d'=3$ and $p_a'=0$. Thus  $|C'|$ is a 4-dimensional linear system whose general member is a smooth rational curve of degree 3. Conversely, if $d=3$ and $p_a=0$ then  $|C'|$ is a 5-dimensional linear system whose general member is a smooth elliptic curve.

(\ref{curve.dp4.say}.3) If $d=4$ and $g=0$ then $d'=4$ and $g'=0$. Thus $|C'|$ is a 3-dimensional linear system whose general member if a smooth rational curve.

(\ref{curve.dp4.say}.4) If $d=6$ and $g=0$ then we formally get that $d'=2$ and $g'=-2$ but
usually there is no such curve $C'$.  We need to understand this case somewhat differently.

Note that  $(C^2)=4$, thus $|C|$ maps $S_4$ onto a degree 4 surface in
$\p^5$. There are 2 such surfaces over an algebraically closed field.

One is the Veronese, that is, $\p^2$ embedded by $|\o_{\p^2}(2)|$. 
This represents $S_4$ as $\p^2$ blown up at 5 points.
Note that in this case twice the conic through the 5 points gives 
the residual curve $C'$. This is a non-generic situation.

The other surface is $\p^1\times \p^1$ embedded by
$|\o_{\p^1\times \p^1}(1,2)|$. This represents $S_4$ as $\p^1\times \p^1$ blown up at 4 points and $|C|$ is the birational transform of 
$|\o_{\p^1\times \p^1}(1,2)|$. 
\end{say}

Working with  degree 4  del~Pezzo surfaces as in 
(\ref{RC.quad.thm}) 
allows us to describe a few more spaces of maps to $Q^3$. 
In all these cases the curves are not linearly normal in $\p^4$,
so these results do not yield   descriptions for higher dimensional quadrics.

\begin{prop}\label{06.16.prop}
 Let $Q^3$ be a smooth quadric. Then
\begin{enumerate} 
\item $\mo_0(6, Q^3)\times \p^3\simb \p^{21}$ and
\item $\mo_1(6, Q^3)\simb \sym^2\bigl(\OG(\p^1, Q^3)\bigr)\times \p^{12}$.
\end{enumerate}
\end{prop}

Proof. The elliptic case is easier here.
Let $E_6\subset Q^3$ be a general, smooth curve of genus 1.
Computing as in the proof of (\ref{RC.quad.thm}) 
we see that $E_6$ is contained in a pencil of quadric sections of $Q^3$
by (\ref{curve.dp4.say}.1). 
Intersecting any 2 of them gives a  pair of lines as the  residual intersection curve.
This gives a map $\mo_1(6, Q^4)\map \sym^2\bigl(\OG(\p^1, Q^3)\bigr)$.

In order to understand the fiber, note that a given pair of disjoint lines
$\ell_1+\ell_2$ is contained in a  $\p^7$-of quadrics hence in a $\p^{12}$
of pencils of quadrics. The residual intersections give degree 6 curves of genus 1. This shows (2).

Next let $C_6\subset Q^3$ be a general, smooth, rational  curve.
As we checked in (\ref{curve.dp4.say}.4), $C_6$ lies on a
unique degree 4  del~Pezzo surface $S$ and on it 
$|C_6|$ is a linear system mapping onto  
a surface $P$ that is geometrically isomorphic to  $\p^1\times \p^1$.
Note, however, that $|C_6|$ corresponds to $|\o_{\p^1\times \p^1}(1,2)|$.
Thus, over $k$, 
$P\cong \p^1\times Q^2$ where $Q^2$ is a conic in some $\p^2$
and $|C_6|$ corresponds to $\pi_1^*\o_{\p^1}(1)\otimes\pi_2^*\o_{\p^2}(1)$ 
where $\pi_i:P\to \p^i$ are the coordinate projections.
Let us now consider 
$$
|\o_{\p^1}(1)|\times |\o_{\p^1}(2)|
\subset |\pi_1^*\o_{\p^1}(1)\otimes\pi_2^*\o_{\p^2}(1)|.
$$
Geometrically, a general member of $ |\o_{\p^1}(1)|\times |\o_{\p^1}(2)|$  consists of 3 lines. On $S_4$ these become
 3 conics
$B_1+B_2+B_3$ where $B_2$ is defined over $k$ and $B_1, B_3$ are disjoint conjugates. So a general point in  $\mo_0(6, Q^3)\times \p^3$
is represented by $(C_6, B_1+B_2+B_3)$.

In order to build this space from the other direction, we first pick
$B_2$, then a pair of points $\{p_1, p_3\}$ on $B_2$
and after that   a pair of conics through $p_1, p_3$.
These choices involve  
$\grass (\p^2, \p^4)$, $\sym^2(Q^2)\cong \p^2$ and then
(the Weil restriction of) $\grass (\p^1, \p^3)$ which parametrizes pairs
of 2-planes containing $\{p_1, p_2\}$. 
Together these give  $\p^{16}$. 
Finally note that $B_1+B_2+B_3$ is a curve of degree 6 and genus 0.
By  (\ref{curve.dp4.say}.3), a generic $B_1+B_2+B_3$ lies on a
unique degree 4  del~Pezzo surface $S$ and on it 
$|B_1+B_2+B_3|$ is a 5-dimensional linear system 
whose general member is a degree 6 smooth rational curve. 
\qed

\medskip

One can try to use K3 surfaces in a similar way. This leads to some
interesting examples but nothing very useful for our present purposes.

\begin{say}[Residual intersections on K3 surfaces]\label{K3.res.say}
It was observed by \cite{Kovacs94} that if a K3 surface contains a smooth rational curve then usually it also contains another one. 
Computing the class of the other curve gives some interesting examples.

  ($C_9\subset Q^3$). 
$H^0\bigl(Q^3, \o_{Q^3}(3)\bigr)$ has dimension $\binom{7}{4}-5=30$.
Thus a general $C_9\subset Q^3$ is contained in a unique pencil of
$| \o_{Q^3}(3)|$.  The intersection of any two members gives
a curve  $C_9+C'_9$ where $C'_9$ is another degree 9 rational curve
meeting $C$ in 29 points.  In fact, $C'_9$  is the unique  degree 9 rational curve in $\p^4$
meeting $C$ in 29 points. Indeed, any such rational curve $C^*_9$
would meet $Q^3$ in $29>2\cdot 9$ points, thus $C^*_9$ is contained in $Q^3$.
Next,  $C^*_9$
would meet any member of the  pencil $| \o_{Q^3}(3)|(-C_9)$
in at least $29>3\cdot 9$ points, thus $C^*_9$ is contained in
the  base locus of  $| \o_{Q^3}(3)|(-C_9)$, which is exactly $C_9\cup C'_9$.

This gives a natural involution on  $\mo_0(Q^3,9)$.

($C_8\subset Q^3$).  Computing as above we see that a general $C_8$
is contained in a $\p^5$ of K3 surfaces of degree 6. Following the method
of \cite{Kovacs94} gives that a general  K3 surface contains another smooth rational curve that is linearly equivalent to $39H-14C_8$. Thus it has degree 200 and meets $C_8$ in 262 points. 
\end{say}

\section{Degenerate maps of $\p^2$}\label{sec.degmaps}

Let $\Phi\in \maps_2(\p^2, Q^n)$ be given by its coordinate functions
$h_0,\dots, h_{n+1}$.
The $h_i$ are in the 6-dimensional vector space $V^6$ of degree 2 homogeneous polynomials in 3 variables. 
We say that $\Phi$ is {\it non-degenerate} if the $h_i$ span $V^6$
and  {\it degenerate} otherwise.  Thus the image of a degenerate map
lies in a (possibly singular) subquadric  $Q^r\subset Q^n$ of dimension $\leq 3$.
In this section we describe degenerate maps $\p^2\map Q^n$, based on their image.

\begin{say}[Expected dimension of $\maps_2(\p^2, Q^n)$]\label{pr.var.say}
A map $\Phi:\p^2\map Q^n$ is a map $\Phi:\p^2\map \p^{n+1}$
whose image lies in $Q^n$. Quadratic maps from $\p^2$ to $\p^{n+1}$
are given $n+2$ quadratic forms
$$
h_0(u,v,w),\dots, h_{n+1}(u,v,w)\in H^0\bigl(\p^2, \o_{\p^2}(2)\bigr).
\eqno{(\ref{pr.var.say}.1)}
$$
We exclude the case when all of the $h_i$ are identically 0
and $(h_0,\dots, h_{n+1})$ is identified with
$(\lambda h_0,\dots, \lambda h_{n+1})$ for any nonzero constant $\lambda$.
Thus these maps naturally form a projective space  of dimension $6(n+2)-1=6n+11$. 
The image lies on $Q^n$ iff
$$
q\bigl(h_0(u,v,w),\dots, h_{n+1}(u,v,w)\bigr)\equiv 0.
\eqno{(\ref{pr.var.say}.2)}
$$
The left hand side is a degree 4 homogeneous polynomial in $u,v,w$
whose coefficients are degree 2 polynomials in the 
coefficients of the $h_i$.
Thus  (\ref{pr.var.say}.2) is equivalent to  $\tbinom{6}{2}=15$ 
quadratic equations in the coefficients of the $h_i$.
Thus we have a natural realization
$$
\maps_2(\p^2, Q^n)\subset  \p^{6n+11}
\eqno{(\ref{pr.var.say}.3)}
$$
as the common zero set of $15$ quadratic equations. In particular, 
$$
\dim \maps_2(\p^2, Q^n)\geq 6n+11-15=6n-4.
\eqno{(\ref{pr.var.say}.4)}
$$
More precisely, every
irreducible component of $\maps_2(\p^2, Q^n) $ has dimension
$\geq 6n-4$. We will see that equality holds if $n\geq 3$. 
\end{say}

\begin{say}[Maps with small dimensional image] \label{conics.say}
Maps with 0-dimensional image are parametrized by $Q^n$.

Maps with 1-dimensional image are obtained as  composites
$\p^2\map \p^1\to Q^n$. The composite can be a quadratic map in two ways.
If $\p^1\to Q^n$ is linear then $\p^2\map \p^1$ should be a pencil of conics,
giving a  moduli space that is birational to 
$\maps_1(\p^1, Q^n)\times \p^8$. 
If $\p^1\to Q^n$ is quadratic then $\p^2\map \p^1$ should be a pencil of lines,
giving a  moduli space that is birational to 
$\maps_2(\p^1, Q^n)\times \p^2$. 
Using Theorem \ref{curve.map.thm},  the spaces of maps whose image 
is a line (resp.\ conic) are  further birational to
$$
\begin{array}{lcl}
\maps_L(\p^2, Q^n)&\simb & \OG(\p^1,Q^n)\times \p^{11} \qtq{and}\\ \maps_C(\p^2, Q^n)&\simb & Q^n\times \p^{2n+2}.
\end{array}
\eqno{(\ref{conics.say}.1)}
$$

 In particular, if $n=1$ and there is a map $\p^2\map Q^1$
then $Q^1(k)\neq \emptyset$, thus we can write
$q\sim x^2-yz$.  
The simplest solution is the inverse of the projection from
$(0,1,0)$. More generally, we get solutions of the form
$$
\begin{array}{rcl}
h_x(u,v,w)&=&\ell_1(u,v,w)\ell_2(u,v,w),\\
h_y(u,v,w)&=&\alpha\ \ell_1^2(u,v,w),\\
h_z(u,v,w)&=&\alpha^{-1}\ell_2^2(u,v,w).
\end{array}
\eqno{(\ref{conics.say}.2)}
$$
\end{say}

\begin{say}[Quadric cone image]\label{q.2-q20.sing.say}
Assume that  the image of $\Phi:\p^2\to Q^n$ is a singular quadric surface.
Let $p\in Q^n$ be the singular point of $\Phi(\p^2)$. Composing with projection from
$p$ we get a map $\pi_p\circ \Phi: \p^2\map Q^{n-2}_p$ whose image is a conic.
By (\ref{conics.say}.1) the later form a space birational to
$Q^{n-2}_p\times \p^{2n-2}$.  Choosing $p\in Q^n$ and $p'\in Q^{n-2}_p$
is the same as choosing the line  $\ell=\langle p, p'\rangle\in \OG(\p^1, Q^n)$
and the point $p\in \ell$.  Thus the maps $\pi_p\circ \Phi$ are
parametrized by $\OG(\p^1, Q^n)\times \p^{2n-1}$.  
Lifting  a map  $\phi^1:\p^2\map Q^1=\bigl(q(x_0,x_1,x_2)=0\bigr)\subset \p^2$ to a 
map to the corresponding quadric cone $\phi^{2,0}:\p^2\map Q^{2,0}=\bigl(q(x_0,x_1,x_2)=0\bigr)\subset \p^3$ 
is equivalent to choosing an arbitrary  3rd component $\phi^{2,0}_3(u,v,w)$.
 Thus  the space of maps with  quadric cone image is   birational to
$$
\maps_{QC}(\p^2, Q^n)\simb \OG(\p^1, Q^n)\times \p^{2n+5}.
\eqno{(\ref{q.2-q20.sing.say}.1)}
$$

\end{say}

\begin{say}[Quadruple planes] \label{quadr.pl.say}
Here we consider the cases when
the  image of $\Phi(\p^2)$ is a linear subspace $L^2\subset Q^n$. Thus, in suitable coordinates
$\Phi$ is given as  $(h_0, h_1, h_2, 0,\dots,0)$. 
These maps are parametrized by 
$h_0, h_1, h_2\in H^0\bigl(\p^2, \o_{\p^2}(2)\bigr)$, up to a multiplicative constant.
For general choices  
$h_0, h_1, h_2 $  have no common zero and then $\Phi$ is an everywhere defined degree 4 morphism
$\p^2\to L^2$. If $n=4$ then  $\Phi_*[\p^2]$ equals  $4[A]$ or $4[B]$
and every deformation again has the same class.
If the $h_i$ have common zeros, we get maps of lower degree
and even maps whose image is a conic, a line or a point.

All  $2$-planes in $Q^n$  are parametrized by the   orthogonal Grassmannian 
 $\OG(\p^2,Q^n)$. 
As we noted before Corollary  \ref{4-folds.3.chow}, the  universal family over 
 $\OG(\p^2,Q^n)$ is birationally trivial.
Thus 
$$
\maps_{QP}(\p^2, Q^n)\simb  \OG(\p^2,Q^n)\times\p^{17}.
\eqno{(\ref{quadr.pl.say}.1)}
$$

If $n=4$ then
there is a 2-plane defined over $k$ iff $q$ is split and in this case
$\OG(\p^2,Q^4)\cong \p^3\amalg \p^3$, hence
$\maps_2(\p^2, Q^4)$ is birational to the disjoint union of 2 copies of $\p^{20}$.
\end{say}

\begin{say}[Smooth quadric surface image]\label{q.2-q20.say}
If $Q^2$ is split then it is isomorphic to $\p^1\times \p^1$ and
quadratic maps $\p^2\map \p^1\times \p^1$ correspond to 
a pair of linear projections  $\p^2\map \p^1$.
If $Q^2=xy-zt$  then these are given by
$$
\begin{array}{lcl}
h_x(u,v,w)&=&\ell_1(u,v,w)\ell_3(u,v,w),\\
h_y(u,v,w)&=&\ell_2(u,v,w)\ell_4(u,v,w),\\
h_z(u,v,w)&=&\ell_1(u,v,w)\ell_4(u,v,w),\\
h_t(u,v,w)&=&\ell_2(u,v,w)\ell_3(u,v,w).
\end{array}
\eqno{(\ref{q.2-q20.say}.1)}
$$
Each of the pairs  $(\ell_1, \ell_2)$ and  $(\ell_3, \ell_4)$
is determined up to a multiplicative scalar, so 
$\maps_2(\p^2, Q^2)\simb \p^{10}$. 

Geometrically,  blow up the 2 points 
$(\ell_1=\ell_2)$ and $(\ell_3=\ell_4)$
 and contract the line $L$ connecting them to a point $p\in Q^2$.
Thus $\Phi^{-1}$ is given by projecting $Q^2$ from $p$ to $\p^2$ and 
 composing it  with an automorphism of $\p^2$.

Thus if the image of $\phi:\p^2\map Q^n$ is a smooth quadric surface
then $\phi$ can be uniquely obtained by composing the inverse of a projection
$\pi_p:Q^n\map \p^n$ with a linear embedding $\p^2\map \p^n$. These 
{\it unprojection maps} form a variety 
$$
\maps_{UP}(\p^2, Q^n)\simb  Q^n\times \p^{3n+2}.
\eqno{(\ref{q.2-q20.say}.2)}
$$ 
In particular, for a quadric surface $Q^2$,  if there is a map $\p^2\map Q^2$
then $Q^2(k)\neq \emptyset$.
This holds iff  $q\sim x^2-ay^2-zt$.  
The inverse of the projection from
$(0,0,1,0)$  gives the obvious solution
$(uw,vw,u^2-av^2, w^2)$. More generally, if we write
$h_x(u,v,w)\pm \sqrt{a}h_y(u,v,w)$ as products of linear factors  
$$
h_x\pm \sqrt{a}h_y=
\bigl(\ell_1\pm \sqrt{a}\ell_2\bigr)
\bigl(\ell_3\pm \sqrt{a}\ell_4\bigr), 
$$
we get the general solutions
$$
h_x=\ell_1\ell_3+a\ell_2\ell_4,\
h_y=\ell_1\ell_4+\ell_2\ell_3,\
h_z=\ell_1^2-a\ell_2^2,\
h_t=\ell_3^2-a\ell_4^2,
\eqno{(\ref{q.2-q20.say}.3)}
$$
but  some non-obvious changes of $\ell_1,\dots, \ell_4$   result in the same map.
\end{say}

Before we go further, we need to recall some
classical facts about the Veronese surface;
see 
\cite[Cap.XV]{bertini-book}, \cite[Kap.XVI]{bertini-book-g} or  
\cite[Sec.VII.3]{MR0034048} for details.

\begin{say}[Veronese surface and its projections]\label{veronese.pr.say}
Let $\Phi:\p^2\map \p^{n+1}$ be a map given by degree 2  homogeneous polynomials. Since the space of degree 2  homogeneous polynomials
in 3 variables has dimension 6, $\Phi(\p^2)$ is always contained in a 5-dimensional linear subspace of $ \p^{n+1}$. 
Moreover, up to $\aut(\p^5)$, there is a unique embedding
$\Phi:\p^2\map \p^5$  given by degree 2  homogeneous polynomials
whose image is not contained in any hyperplane. Its 
image is called a {\it Veronese surface;} we denote it by 
 $V\subset \p^5$. A typical example is
$$
\Phi: (u:v:w)\mapsto  (u^2: v^2: w^2:uv:vw:uw).
$$
For every other quadratic map $\Phi':\p^2\map \p^{n+1}$,
we can think of $\Phi'(\p^2)$  as a projection of the
Veronese surface.

If we think of the dual $(\p^5)^{\vee}$ as the space of conics in $p^2$,
then we see that 
$\aut(\p^2)$ acts on $\p^5$ with 3 orbits, corresponding to the rank if the conic. 
The rank 1 orbit is the Veronese surface, the rank 2 orbit is its secant variety  and the rank 3 orbit is open.

This makes it easy to determine the projections of $V$ to $\p^4$: we just need to compute one example for each of the 3 orbits.

$\bullet$ If we project from a point $p$ not on the secant variety,
the projection is an isomorphism and we get a smooth degree 4 surface
$V_1$. For example, projecting from $(1:1:1:0:0:0:0)$ gives
$$
(u:v:w)\mapsto  (u^2-w^2:  v^2-w^2:uv:vw:uw).
$$
We see that $V_1$  does not lie on any quadric 3-fold $Q^3$.

$\bullet$ If we project from a point  on the secant variety but not on $V$,
the projection is a degree 4 surface with a double line 
$V_2$. In suitable coordinates
it is the image of the map  
$$
(u:v:w)\mapsto  (u^2: v^2: w^2:uv:vw).
$$
It is contained in the family of quadrics
$\lambda (x_0x_1-x_4^2)+\mu(x_1x_2-x_5^2)$. All of these are singular,
namely at the  point $(\mu:0:-\lambda:0:0)$.

Thus $V_2$ is a complete intersection of any 2 quadrics that contain it.

$\bullet$ If we project from a point  on $V$,
the projection is a degree 3 surface 
$V_3$.  In suitable coordinates
it is the image of the map  $\p^2\map \p^4$ given by
$$
(u:v:w)\mapsto  (u^2: v^2:uv:vw:uw).
$$
The map blows up the point  $(0:0:1)$. The image
 is a degree 3 ruled surface.
$V_3$ lies on a 3-dimensional family of quadrics best described 
as the $2\times 2$ subdeterminants of the matrix
$$
\left(
\begin{array}{ccc}
x_0 & x_2 & x_4\\
x_2 & x_1 & x_3
\end{array}
\right).
$$
\end{say}

\begin{say}[Quadric 3-folds]\label{3f.say} We have seen on our list 
in (\ref{veronese.pr.say}) that a
smooth quadric 3-fold can not contain a degree 3 or 4 projection of a Veronese surface. Thus for all non-degenerate maps the image
is a quadric surface contained in $Q^3$.

Maps whose image is a  smooth $Q^2\subset Q^3$ were described in 
(\ref{q.2-q20.say}.3) and 
$$
\maps_{UP}(\p^2, Q^3)\simb  Q^3\times \p^{11}.
\eqno{(\ref{3f.say}.1)}
$$ 
 This family has members over $k$ whenever
$Q^3(k)\neq \emptyset$ and then
$\maps_{UP}(\p^2, Q^3)\simb \p^{14}$.

Maps to a quadric cone
$Q^{2,0}\subset Q^3$ form another family of dimension
$3+11=14$.  By (\ref{q.2-q20.sing.say}.1) they form a family
$$
\maps_{QC}(\p^2, Q^3)\simb \OG(\p^1, Q^3)\times \p^{11}.
\eqno{(\ref{3f.say}.2)}
$$
The second family has members over $k$ whenever
$Q^3$ contains a $k$-line. This happens only if
$$
q(x_0,\dots, x_4)\sim y_0^2+y_1y_2+y_3y_4.
$$
Both of these families have dimension 14, so we conclude that
$$
\maps_2(\p^2, Q^3)\simb \bigl(Q^3\times \p^{11}\bigr)\amalg  \bigl(\OG(\p^1, Q^3)\times \p^{11}\bigr).
\eqno{(\ref{3f.say}.3)}
$$
\end{say}

\begin{say}[Projected Veronese surfaces]\label{pvs.map.say}
Let $V_2\subset Q^n$ be a projected Veronese surface
with double line $\ell\subset V_2$.

 The linear span $\langle V_2\rangle$ of $V_2$ has
dimension 4 and $V_2$ lies in the 3-dimensional quadric
$Q^n\cap \langle V_2\rangle$. As we noted in (\ref{3f.say}),  $Q^n\cap \langle V_2\rangle$ has a unique  singular point $p$.
Thus $V_2$ determines a singular subquadric
  $Q^{3,0}\subset Q^n$ plus a line $\ell\subset Q^{3,0}$ through the vertex $p$
of $Q^{3,0}$.
We can specify these data by first picking the  line  $\ell$,
then the point    $p\in \ell$ and finally a 3-dimensional linear subspace
$[\ell]\in L^3\subset \p(T_pQ^n)$.

Fix  $Q^{3,0}:=(x_1x_2=x_3x_4)\subset \p^4$ and the line $\ell=(x_2=x_3=x_4=0)$. 
One easily computes that quadric sections of $Q^{3,0}$ that vanish doubly along $\ell$ are all of the form $x_2L(x_0,\dots, x_4)+a_3x_3^2+ a_4x_4^2$, thus their 
linear system has dimension 6.

Putting these together we obtain that the Chow variety of
 projected Veronese surfaces is
$$
\begin{array}{rcl}
\chow_{PV}(Q^n)&\simb& \OG(\p^1, Q^n)\times \p^1\times \p^{3(n-4)}\times \p^6\\
&\simb& \OG(\p^1, Q^n)\times \p^{3n-5}.
\end{array}
\eqno{(\ref{pvs.map.say}.1)}
$$
Since the quadratic maps $\p^2\to V_2$ are birational,
we conclude that
$$
\maps_{PV}(\p^2, Q^n)\simb \OG(\p^1, Q^n)\times \p^{3n+3}.
\eqno{(\ref{pvs.map.say}.2)}
$$ 
\end{say}

\section{Veronese surfaces in quadric 4-folds}\label{sec.veron}

In the previous section we have  enumerated all degenerate maps of $\p^2$
to a smooth quadric 4-fold.
Among those, $\maps_{PV}(\p^2, Q^4)$ and  $\maps_{QP}(\p^2, Q^4)$
have dimension $20$ and the others have dimension $<20$. 
Since the expected dimension of $\maps_2(\p^2, Q^4)$ is $6\cdot 4-4=20$ by
(\ref{pr.var.say}), once we prove that  $\dim\maps_{V}(\p^2, Q^4)=20$,
we will know that all the geometric irreducible components of
$\maps_2(\p^2, Q^4)$ are given by
 $\maps_{V}(\p^2, Q^4), \maps_{PV}(\p^2, Q^4)$ and  $\maps_{QP}(\p^2, Q^4)$.

 Thus it remains to describe the moduli space of Veronese surfaces in quadric 4-folds.

\begin{say}[Homology class of Veronese surfaces]
The Chern class of  $Q^4$ is  $(1+H_5)^6(1+2 H_5)^{-1}$ where $ H_5$ is the hyperplane class on $\p^5$. Restricting to
$V$ this gives  $(1+2H_2)^6(1+4H_2)^{-1}$ where $ H_2$ is the hyperplane class on $\p^2$ (note that the restriction of $ H_5$ equals $2H_2$).
The  Chern class of  $\p^2$ is  $(1+H_2)^3$, hence the Chern class of the normal bundle of $V\subset Q^4$ is 
$$
c\bigl(N_{V,Q}\bigr)=(1+2H_2)^6(1+4H_2)^{-1}(1+H_2)^{-3}=1+5H_2+10H_2^2.
$$
Thus  $(V\cdot V)=10$.  Since $(V\cdot  H_5^2)=4$
this leads to the possibilities  $[V]=3[A]+[B]$ or $[V]=[A]+3[B]$,
as computed in \cite[Sec.XV.13]{bertini-book}, \cite[Sec.XVI.13]{bertini-book-g}.
This suggests, and we will see,  that there are 2 distinct Veronese families. 
\end{say}

We start with a duality between the 2 families of Veronese surfaces.
Then we work with a special case when one side is a very singular
specialization of Veronese surfaces to get our moduli description.

\begin{say}[Linking the two Veronese families]\label{vero.link.say}
Let  $Q^4$  be a smooth quadric 4-fold and $V\subset Q^4$ a
Veronese surface.   $V$ is contained in a 5-dimensional family
of quadrics in $\p^5$, thus in a 4-dimensional family
of quadric sections of $Q^4$. Let
$\langle Q_1, Q_2\rangle \subset |\o_{Q^4}(2)(-V)|$ be
any pencil. Then $Q^4\cap Q_1\cap Q_2$ is a  degree 8 surface containing $V$.
Thus we get a residual intersection
$$
 Q^4\cap Q_1\cap Q_2=V\cup V',
\eqno{(\ref{vero.link.say}.1)}
$$
 where $V'$ is another degree 4 surface contained in $Q^4$.
If  $[V]=3[A]+[B]$ then $[V']=[A]+3[B]$, so the 
  correspondence  $V\leftrightarrow V'$ interchanges the 2 families. 
Most likely, a given $V$ will give all others by repeating such linkages but I have no proof of this.

Let us see  some nice examples, using  the representation of $V$ in the space of symmetric matrices. Our first pair is
$$
V:=\left(
\begin{array}{ccc}
u^2 & uv & uw \\
uv & v^2 & vw \\
uw & vw & w^2
\end{array}
\right)
\qtq{and}
V':=\left(
\begin{array}{ccc}
u^2 & -uv & -uw \\
-uv & v^2 & -vw \\
-uw & -vw & w^2
\end{array}
\right).
\eqno{(\ref{vero.link.say}.2)}
$$
Their union is defined by  the symmetric subdeterminant  equations
$$
x_{01}^2-x_{00}x_{11}=
x_{02}^2-x_{00}x_{22}=
x_{12}^2-x_{11}x_{22}=0.
\eqno{(\ref{vero.link.say}.3)}
$$
The sum of the 3 equations gives
$$
Q^4:=\bigl(x_{01}^2+x_{02}^2+x_{12}^2-x_{00}x_{11}-x_{00}x_{22}-x_{11}x_{22}=0\bigr);
\eqno{(\ref{vero.link.say}.4)}
$$
it is an ellipsoid over $\r$; see  (\ref{sur.circ.all.thm}.4--5). 

We can also take the 3 non-symmetric subdeterminant  equations
$$
x_{01}x_{12}-x_{02}x_{11}=
x_{12}x_{02}-x_{01}x_{22}=
x_{01}x_{22}-x_{12}x_{02}=0.
\eqno{(\ref{vero.link.say}.5)}
$$
In this case $V''$ is the union of four 2-planes, consisting of all matrices of the form
$$
\left(
\begin{array}{ccc}
* & 0 & 0 \\
0 & * & 0 \\
0 & 0 & *
\end{array}
\right),\
\left(
\begin{array}{ccc}
* & * & 0 \\
* & * & 0 \\
0 & 0 & 0
\end{array}
\right),\
\left(
\begin{array}{ccc}
0 & 0 & 0 \\
0 & * & * \\
0 & * & *
\end{array}
\right)\ \mbox{and}\
\left(
\begin{array}{ccc}
* & 0 & * \\
0 & 0 & 0 \\
* & 0 & *
\end{array}
\right).
\eqno{(\ref{vero.link.say}.6)}
$$

Note that the canonical class of the intersection of 3 quadrics in $\p^5$ is trivial, so $V\cup V'$  is a   singular K3 surface.
Another way to construct these  $V\cup V'$ is the following.

A degree 6 elliptic normal curve $\tau:E\into \p^5$
is contained in 4 different Veronese surfaces, an observation going back to Coble \cite{MR1501210}. Indeed, for a given elliptic  curve $E$
its 
degree 6 embeddings $\tau_6:E\into \p^5$ correspond to degree 6 line bundles
$L_6:=\tau_6^*\o_{\p^5}(1)$ and degree 3 embeddings $\tau_3:E\into \p^2$ correspond to degree 3 line bundles
$L_3:=\tau_3^*\o_{\p^2}(1)$. Combining $\tau_3$ with the Veronese embedding gives
a degree 6 elliptic normal curve such that $L_6=L_3^{\otimes 2}$.
Given $L_6$, there are 4 different $L_3$ satisfying this equation. These give
the 4 different Veronese surfaces containing $\tau_6(E)$. 
The union of any 2 of them is a singular K3 surface.

Russo pointed out that the various degenerations of Veronsese surfaces can best be seen using the dual picture to be discussed in (\ref{vero.eq.say}). 

\end{say}

Next we study a special case when $V'$ is singular, even reducible, but
slightly better than the ones exhibited in (\ref{vero.link.say}.5--6).

\begin{say}[A degenerate version]\label{Ver.final.pf.say}
Let $p\in V\subset Q^4\subset \p^5$ be a pointed Veronese surface.
Instead of working with all quadrics that contain $V$, we consider only
quadric cones with vertex $p$ that contain $V$.
These  form a linear system  $|2H_{5}|(-V-2p)$ of dimension 2.  We claim that if $\langle Q_1, Q_2\rangle\subset  |2H_{5}|(-V-2p)$ is a general pencil then
$$
Q^4\cap Q_1\cap Q_2=V\cup P_1\cup P_2\cup Q^2,
\eqno{(\ref{Ver.final.pf.say}.1)}
$$
where  the $P_i$ are 2-planes meeting at $p$ and $Q^2$ is a smooth quadric surface through $p$. This is easier to see by projecting everything from $p$.
So we have $\pi:\p^5\map \p^4$ and 
 $\pi(V)$ is a degree 3 surface $F_3\cong \f_1$, as discussed in (\ref{veronese.pr.say}). 
Thus $|2H_{4}|(- F_3)$ has dimension 2 and if 
$\langle Q'_1, Q'_2\rangle\subset  |2H_{4}|(- F_3)$ is a general pencil then
$$
Q'_1\cap Q'_2=F_3\cup L^2,
\eqno{(\ref{Ver.final.pf.say}.2)}
$$
where $L^2$ is a 2-plane that meets $F_3$ in a conic. This $L^2$ is the projection of the $Q^2$. 

The planes $P_i$ are less visible after projection; we obtain them as follows.
Write $(p\in V)$ as the image of $\phi:(p_0,\p^2)\into \p^5$.
Note that $\phi^{-1}\bigl(V\cap T_pQ^4\bigr)$ is a degree 2 curve in $\p^2$ that is singular at $p_0$. Thus $\phi^{-1}\bigl(V\cap T_pQ^4\bigr)$ is a 
pair of lines  in $\p^2$ that meet at $p_0$ and
$V\cap T_pQ^4$ is a pair of conics $C_i$ meeting at $p$.
Each  $C_i$  spans a 2-plane $P_i$ which is contained every member of $|2H_{\p^5}|(-V-2p)$.
Thus the $P_i$  do not depend on $\langle Q_1, Q_2\rangle$ but $Q^2$ does. 

Note that the $\pi$-images of the $P_i$ are lines $L_i\subset F_3$ and 
$L^2$ meets both of them. Thus $P_1\cup P_2\cup Q^2$ is a surface where
all 3 components meet at $p$,  $P_i\cap Q^2$ are lines and
$P_i\cap P_2=\{p\}$.
These data specify a unique  isomorphism class of surfaces in $\p^5$, up to $\aut(\p^5)$.

This construction can be reversed. 
That is, assume that we have $P_1\cup P_2\cup Q^2\subset Q^4$.
Let $|2H_{\p^5}|(- Q^2-P_1- P_2-2p)$ denote that linear system of 
quadric cones with vertex $p$ that contain $P_1\cup P_2\cup Q^2$.
We see that it has dimension $2$. (This is best seen after projecting from $p$.
Then we have the linear system of quadrics that contain a 2-plane $L^2$ and two lines $L_i$ that meet it.)

If $\langle Q_1, Q_2\rangle\subset |2H_{5}|(- Q^2-P_1- P_2-2p)$
is a general pencil then 
$$
Q^4\cap Q_1\cap Q_2=V\cup P_1\cup P_2\cup Q^2,
\eqno{(\ref{Ver.final.pf.say}.3)}
$$
where $V$ is a Veronese surface passing through $p$. 

Let us next study  the induced correspondence
$$
\{V\}\leftarrow  \{V\cup Q^2\cup P_1\cup P_2\}\rightarrow \{P_1\cup P_2\}.
\eqno{(\ref{Ver.final.pf.say}.4)}
$$ 

We start with $[V]\in \chow_V(Q^4)$. As we noted before Corollary  \ref{4-folds.3.chow}, the  universal family over 
$\chow_V(Q^4)$ is birationally trivial, thus pointed
Veronese surfaces are birationally parametrized by $\chow_V(Q^4)\times \p^2$.
The choice of the pencil $\langle Q_1, Q_2\rangle\subset  |2H_{\p^5}|(-V-2p)$
is another $\p^2$-factor. Thus the left hand side of (\ref{Ver.final.pf.say}.4) is birational to
$$
\chow_V(Q^4)\leftarrow  \chow_V(Q^4) \times \p^4.
$$ 
On the right hand  side of (\ref{Ver.final.pf.say}.4) we start by choosing $P_1\cup P_2 $.
This is the same as picking $P_1, P_2\subset Q^4$ that are in the same geometric irreducible component 
of $\OG(\p^2, Q^4)$. This will be denoted by
$\sym^{2}_{a}\bigl( \OG(\p^2, Q^{4})\bigr)$ in
(\ref{sq.og.some.say}). 
Then we choose $Q^2$. As we noted above, this is equivalent to
having 2 lines  $L_1, L_2\subset \p^4$ and choosing a 2-plane $L^2$ that meet both of them. After fixing an auxiliary $\p^2_2\subset \p^4$, these are
parametrized by $L_1\times L_2\times \p^2$. Note, however, that for us
the $L_1, L_2$ are conjugate, thus instead of  $L_1\times L_2$ we pick 
another auxiliary $\p^2_0\subset \langle L_1, L_2\rangle\cong \p^3$
and represent a point pair $\{p_1, p_2: p_i\in L_i\}$
by the intersection point $\langle p_1, p_2\rangle\cap \p^2_0$. 
(This is just a birational construction of the Weil restriction 
$\Re(L_1)$.)
As before,  the choice of the pencil $\langle Q_1, Q_2\rangle\subset  
|2H_{5}|(- Q^2-P_1- P_2-2p)$
is a $\grass(\p^1, \p^4)$-factor. 
Thus the diagram (\ref{Ver.final.pf.say}.4) is birational to
$$
\begin{array}{ccc}
\chow_V(Q^4) \times \p^4
&\simb& \sym^{2}_{a}\bigl( \OG(\p^2, Q^{4})\bigr) \times \p^{10}\\
\downarrow && \downarrow\\
\chow_V(Q^4) &&
 \sym^{2}_{a}\bigl( \OG(\p^2, Q^{4})\bigr).
\end{array}
\eqno{(\ref{Ver.final.pf.say}.5)}
$$ 
We finally use   (\ref{sq.og.some.say}.4) to get that 
$$
\chow_V(Q^4)\times \p^4\simb 
 (t^2+\Delta=0)\times Q^4\times \p^{12}.
\eqno{(\ref{Ver.final.pf.say}.6)}
$$
\end{say}

\begin{rem}\label{other.comp.rem}
Another family of degree 4 surfaces in $Q^4$ can be obtained  starting with 2 disjoint planes  $P_1, P_2\subset Q^4$.
 A general member of this family is 
$\p^1\times \p^1$ embedded by $\o_{\p^1\times \p^1}(2,1)$. It is better to view such a surface  as the product of a conic with $\p^1$, thus we have
$$
C\times \p^1\subset\p^2\times \p^1\stackrel{j}{\hookrightarrow} \p^5,
$$
where $j$ is  the embedding given by $\o_{\p^2\times \p^1}(1,1)$.
The residual intersection of $Q^4$ with  $j\bigl(\p^2\times \p^1\bigr)$ is a 
pair of planes  $\p^2\times \{p_1,p_2\}$. Conversely, given a 
disjoint pair of conjugate planes  $P_1, P_2\subset Q^4$ and an isomorphism
$\tau:P_1\cong P_2$, we get an embbeding 
$P_1\amalg P_2\subset j\bigl(\p^2\times \p^1\bigr)$.
The  residual intersection of $j\bigl(\p^2\times \p^1\bigr)$ with $Q^4$ gives  a required degree 4 surface
in $Q^4$. This shows that this irreducible component of $\chow(Q^4) $
is birational to  $\sym^{2}_{s}\bigl( \OG(\p^2, Q^{4})\bigr)  \times \p^8$
where 
$\sym^{2}_{s}\bigl( \OG(\p^2, Q^{4})\bigr)$ is the irreducible component of
$\sym^{2}\bigl( \OG(\p^2, Q^{4})\bigr)$ parametrizing disjoint pairs of 2-panes; see (\ref{sq.og.some.say}) for details.

\end{rem}

\section{Quadratic maps up to coordinate changes}\label{sec.coo.ch}

Here we study $\maps_2(\p^2, Q^4)$, up to coordinate changes by 
$\aut(\p^2)\times\aut(Q^4)$. 
The 3 types of irreducible components of $\maps_2(\p^2, Q^4)$
behave differently. There is especially nice geometry behind the Veronese maps.
The following classical facts can be found in
\cite[Cap.XV]{bertini-book}, \cite[Kap.XVI]{bertini-book-g},
 \cite[p.188]{MR0034048} or \cite{ein-she}.

\begin{say}[Equations of Veronese surfaces]\label{vero.eq.say}
The equations of the  Veronese surface are especially clear if we think of the ambient $\p^5$ as the projectivization of the space of $3\times 3$ symmetric matrices. Thus
we consider  $\p^2_{\mathbf u}$ with coordinates  $\{u_0,u_1,u_2\}$ 
and $\p^5_{\mathbf x}$ with coordinates  $\{x_{ij}\}$  for $0\leq i\leq j\leq 2$. 
The Veronese surface $V=V_{\mathbf x}\subset \p^5_{\mathbf x}$ is the image of the map
$\phi:\p^2_{\mathbf u}\into \p^5_{\mathbf x}$ given by
$$
x_{ij}=u_iu_j  \qtq{for} 0\leq i\leq j\leq 2. 
\eqno{(\ref{vero.eq.say}.1)}
$$
In order to describe the equations defining $V$, consider the matrix
$$
M:=
\left(
\begin{array}{ccc}
x_{00} & x_{01} & x_{02} \\
x_{01} & x_{11}& x_{12} \\
x_{02} & x_{12}& x_{22}
\end{array}
\right).
\eqno{(\ref{vero.eq.say}.2)}
$$
$V$ is defined by the equations
$$
\rank M\leq 1,\qtq{equivalently,}  x_{ij}x_{k\ell}=x_{i\ell}x_{kj},
\eqno{(\ref{vero.eq.say}.3)}
$$ 
(where we set $x_{ij}=x_{ji}$)
and  the secant variety of $V$ is  given by
$\det M=0$.
The quadrics containing $V$ are given by (linear combinations of) the $2\times 2$ subdeterminants of $M$. 
These define a birational map
$$
\sigma_{\mathbf x}: \p^5_{\mathbf x}\map \p^5_{\mathbf y}
\qtq{given by} 
y_{ij}=(-1)^{i+j} M_{ji}.
\eqno{(\ref{vero.eq.say}.4)}
$$
Geometrically, we obtain  $\sigma_{\mathbf x} $ by blowing up
 $V_{\mathbf x}\subset \p^5_{\mathbf x}$ and then contracting the 
(birational transform of the) secant variety of 
$V_{\mathbf x}\subset \p^5_{\mathbf x}$.
Note that $\sigma_{\mathbf x}: \p^5_{\mathbf x}\map \p^5_{\mathbf y}$ and
$\sigma_{\mathbf y}: \p^5_{\mathbf y}\map \p^5_{\mathbf x}$
are inverses of each other since,  for an $r\times r$ matrix, 
$$
\operatorname{adj}\big(\operatorname{adj}(M))=(\det M)^{r-2}\cdot M.
\eqno{(\ref{vero.eq.say}.5)}
$$
(For this reason, one usually views $\sigma:=\sigma_{\mathbf x}$
as an involution of $\p^5$, but for us it will be convenient to 
distinguish the source $\p^5$ from the target $\p^5$.)
Note that  $\sigma_{\mathbf x}$ maps quadrics that contain $V$ to hyperplanes;
this correspondence is given as follows. 
Restricting the above factorization of $\sigma_{\mathbf x}$ to a quadric $Q^4$
we obtain the following.  First  we blow up $V_{\mathbf x}\subset Q^4$ and then contract the 
(birational transform of the) secant variety of $V_{\mathbf x}$ in $Q^4$; that is,
the union of those secant lines of $V_{\mathbf x}$ that are contained in $Q^4$.
For the inverse, we blow up $V_{\mathbf y}\cap H^4$ and then contract the 
(birational transform of the) secant variety of $V_{\mathbf y}\cap H^4$.
If $Q^4$ is smooth then the blow-up of $V_{\mathbf x}$ is smooth
but the blow-up of $V_{\mathbf y}\cap H^4$ is smooth iff
$V_{\mathbf y}\cap H^4$ is smooth.  Thus we get the following.
\medskip

{\it Claim \ref{vero.eq.say}.6.} 
  $\sigma_{\mathbf x}$ establishes a one-to-one correspondence between
\begin{enumerate}
\item[(a)] pairs $V_{\mathbf x}\subset Q^4$, up to isomorphism and
\item[(b)] pairs $V_{\mathbf y}\cap H^4\subset  H^4$,  up to isomorphism.
\end{enumerate}
Furthermore, $Q^4$ is smooth iff the corresponding $V_{\mathbf y}\cap H^4$ is smooth. \qed
\medskip

I have not been able to track down the original author of this claim.
For the geometric version
all the ingredients are in \cite[Secs.XV.9--10]{bertini-book},
 \cite[Thm.17]{MR1501210}
 and   
 it is  mentioned as an exercise  in 
 \cite[p.188]{MR0034048}.

Next observe that $V_{\mathbf y}\cap H^4$ is a degree 4 geometrically rational curve in $\p^4$, thus it is obtained as the image of a plane conic  $C$
embedded by $|-2K_C|$. In particular, the pair $V_{\mathbf y}\cap H^4\subset  H^4$
is uniquely determined by 
$V_{\mathbf y}\cap H^4$,  up to isomorphism.
This proves Complement \ref{V.UP.TO.ISOM.compl}. 
\end{say}

One can write down the correspondence  between conics and
Veronese surfaces contained in a smooth quadric as follows.

\begin{say}[Explicit formulas]\label{vero.eq.ex.say} 
Any plane conic
can be given by a diagonal equation
$$
C:=\bigl(a v_0^2+b v_1^2+ c v_2^2=0\bigr)\subset \p^2_{\mathbf v}.
\eqno{(\ref{vero.eq.ex.say}.1)}
$$
This conic is cut out of the Veronese surface $\p^2_{\mathbf v}\into \p^5_{\mathbf y}$,
by the hyperplane
$$
H^4_C:=\bigl(a y_{00}+b y_{11}+ c y_{22}=0\bigr)\subset \p^5_{\mathbf y},
\eqno{(\ref{vero.eq.ex.say}.2)}
$$
and taking its  preimage by $\sigma_{\mathbf x}$ gives the quadric
$$
Q^4_C:=\bigl(a ( x_{11}x_{22}-x_{12}^2)+  
b (x_{00}x_{22}-x_{02}^2)+  
c (x_{00}x_{11}-x_{01}^2) =0\bigr)\subset \p^5_{\mathbf x}.
\eqno{(\ref{vero.eq.ex.say}.3)}
$$
We can assume that $c=1$ and  rearrange the equation into
$$
\bigl(x_{00}+ax_{22}\bigr)
\bigl(x_{11}+bx_{22}\bigr) = x_{01}^2+
ax_{12}^2+bx_{02}^2+abx_{22}^2.
\eqno{(\ref{vero.eq.ex.say}.4)}
$$
Thus we see that  $Q^4_C$ has a $k$-point
and its Witt reduction is the quadric surface
$$
Q^2_C:=\bigl(x_{01}^2+
ax_{12}^2+bx_{02}^2+abx_{22}^2=0)\cong C\times C,
\eqno{(\ref{vero.eq.ex.say}.5)}
$$
where the isomorphism holds by Paragraph \ref{prod.of.conics}.
This shows that the conic $C$ 
and the quadric $Q^4_C$   uniquely determine each other, up to isomorphism. 

We can also use the  the symmetric subdeterminant  equations
(\ref{vero.link.say}.3) to get the quadric
$$
Q^4_{a,b}:=a\bigl(x_{01}^2-x_{00}x_{11}\bigr)+
b\bigl(x_{02}^2-x_{00}x_{22}\bigr)+
ab\bigl(x_{12}^2-x_{11}x_{22}\bigr)=0.
\eqno{(\ref{vero.eq.ex.say}.6)}
$$
The equation can be rewritten in the usual form as 
$$
x_{00}^2+ax_{01}^2+bx_{02}^2+abx_{12}^2=(x_{00}+ax_{11})(x_{00}+bx_{22}).
\eqno{(\ref{vero.eq.ex.say}.7)}
$$

\end{say}

\begin{say}[Decomposable quadric surfaces]\label{prod.of.conics}  Let $C:= (ax_0^2+bx_1^2+cx_2^2=0)$.  Then
$$
\begin{array}{l}
\psi: (x_0:x_1:x_2)\times (y_0:y_1:y_2)\\
\qquad \mapsto
\bigl( x_1y_2-x_2y_1: -x_0y_2+x_2y_0: x_0y_1-x_1y_0: ax_0y_0+bx_1y_1+cx_2y_2
\bigr)
\end{array}
$$
is an isomorphism of $C\times C$ and  the quadric
$$
Q:=(bcz_0^2+acz_1^2+abz_2^2+z_3^2=0).
$$
Note that $\psi$ is given by  4  sections of $\o_{C\times C}(2,2)$ that vanish along the diagonal.

Explanation. The first 3 components of $\psi$  map 2 points on $C$ to the line connecting them. This gives a morphism to (the dual)  $\p^2$ and its branch locus is the set of tangent lines of $C$. Tangent lines have coefficients  $(2ax_0,2bx_1, 2cx_2)$ and these satisfy the equation $bc(2ax_0)^2+ac(2bx_1)^2+ab(2cx_2)^2=0$. This shows that
$bcz_0^2+acz_1^2+abz_2^2$ is (a constant times) a square. 

The inverse of $\psi$  is  given by
$$
\begin{array}{lcl}
(x_0:x_1:x_2) & = &  \bigl(z_1z_3-bz_0z_2: -z_0z_3-az_1z_2: bz_0^2+az_1^2\bigr)\qtq{and}\\[1ex]
(y_0:y_1:y_2) & = &  \bigl(z_1z_3+bz_0z_2: +z_0z_3-az_1z_2: bz_0^2+az_1^2\bigr).
\end{array}
$$
\end{say}

Next we prove the other  claims about maps up to coordinate changes.

\begin{say}[Proof of Proposition \ref{QP.UP.TO.ISOM.THM}]   \label{QP.UP.TO.ISOM.pf}
If $Q^4\subset \p^5$ contains a  plane, we may assume that it is 
$(x_3=x_4=x_5=0)$. Then  $\p^2\to Q^4$ is of the form 
$$
(u:v:w)\mapsto 
\bigl(q_0(u,v,w), q_1(u,v,w), q_2(u,v,w),0,0,0\bigr).
$$
where the $q_i$ are homogeneous of degree 2.
The $q_i$ span a  base point free  net of conics and, if   $\chr k\neq 3$ 
then any such  net  can be written as
$$
\Bigl\langle\tfrac{\partial C(u,v,w)}{\partial u}, \tfrac{\partial C(u,v,w)}{\partial v}, \tfrac{\partial C(u,v,w)}{\partial w}\Bigr\rangle,
$$
where $C(u,v,w)$ is a homogeneous cubic;
 see \cite[p.42]{MR0034048}.

 Thus the  equation of $Q^4$ can be written as
$x_0\ell_3+x_1\ell_4+x_2\ell_5=0$ where the $\ell_i$ are linear. 
Then $Q$ is smooth iff the 6 linear terms in the equation are
linearly independent, hence we can change $x_3, x_4, x_5$ to   $x_i:=\ell_i$ for $i=3,4,5$.

This gives an isomorphism between the  moduli space of quadruple planes
(up to  $\aut(\p^2)\times\aut(Q^4)$)
and  the moduli space of
 plane cubics (up to $\aut(\p^2)$). 
This completes the proof of Proposition \ref{QP.UP.TO.ISOM.THM}. 

Note that we have in fact proved the much more precise statement:
there is an isomorphism of stacks
$$
\aut(\p^2)\big\backslash\maps_{QP}^{\circ}(\p^2, Q^4)\big/\aut(Q^4)
\cong
\Div^{\circ}_{3}(\p^2)\big/\aut(\p^2)
$$
where $\Div^{\circ}_{3}(\p^2)$ parametrizes smooth, degree 3 plane curves.
Equivalently, for any field $k$ we have an equality of sets
$$
\PGL_3(k)\big\backslash\maps_{QP}^{\circ}(\p^2, Q^4)(k)\big/\OO(Q^4,k)
\cong
\Div^{\circ}_{3}(\p^2)(k)\big/\PGL_3(k). 
$$
\end{say}

\begin{say}[Proof of Proposition \ref{PV.UP.TO.ISOM.THM}] Let $S\subset Q^4$ be a projected Veronese surface with singular line $L\subset S$.
Then $S$ is contained in a unique $Q^{3,0}\subset Q^4$ and projecting $S$ from the vertex of  $Q^{3,0}$ gives a birational map  $\pi:S\to Q^2$. 
The image of the singular line in $S$ is a point  $p\in Q^2$.

Conversely, given  $p\in Q^2$ we first  blow up $p$ 
to get $B_pQ^2$ with exceptional curve $E$ 
and then 
obtain $\tau:B_pQ^2\to S$  that is an isomorphism 
$B_pQ^2\setminus E\cong S\setminus L$ outside $E$ 
but whose restriction $\tau_E:E\to L$ has degree 2. We can identify
$\tau_E:E\to L$ with a degree 2 homogeneous polynomial on
$T_pQ^2$, up to scalars.

Note that the quadric $Q^2$ itself gives a quadratic form on $T_pQ^2$;
in suitable coordinates it is $x^2+\Delta y^2$ where $\Delta$ is the discriminant of $Q^4$. Thus, up to isomorphisms, $S$ is determined by
$\sym^2(\p^1)/\OO(x^2+\Delta y^2)$.
 \qed 

\end{say}

\begin{say}[Proof of Proposition \ref{V.UP.TO.ISOM.5.THM}]
\label{pf.V.UP.TO.ISOM.5.THM}
Assume that $Q^5$ contains a Veronese surface $V$. The linear span of $V$ intersects $Q^5$ in a 4-dimensional subquadric $Q^4_V$. As we noted in Paragraph \ref{veronese.pr.say}, a Veronese surface is not contained in any 4-dimensional quadric with an isolated singularity, thus  $Q^4_V$ is smooth.

Fix now an equation $q(x_0,\dots, x_6)=0$ for $Q^5$. 
Using Proposition \ref{PV.UP.TO.ISOM.THM},
after a  coordinate change we can write $q$ as
$$
A\bigl(x_0x_5+x_1^2+ax_2^2+bx_3^2+abx_4^2\bigr)+Cx_6^2
$$
where $x_6=0$ defines $Q^4_V$. Note that the discriminant of the
form in parenthesis is $-1$, thus $-C/\Delta(q)$ is a square.
We can thus choose $C=-\Delta(q)$, independent of the original choice of $V\subset Q^5$. By Witt's cancellation theorem  (\ref{witt.cancel})
this shows that $Q^4_V$ is uniquely determined by $Q^5$, up to isomorphism.
Thus  Proposition \ref{PV.UP.TO.ISOM.THM} shows that  the pair
$(V\subset Q^5)$ is unique  up to isomorphism.\qed
\end{say}

\section{Quadrics of dimension $\geq 5$}\label{sec.5.dim}

In Section \ref{sec.degmaps} we proved that every irreducible component of
$\maps_2(\p^2, Q^n)$ has dimension $\geq 6n-4$ and we computed the spaces that parametrize degenerate maps. In all cases they have dimension $<6n-4$ provided $n\geq 5$. Thus $\maps_V^{\circ}(\p^2, Q^n)$ is dense in $\maps_2(\p^2, Q^n)$
 for every $n\geq 5$. 
In particular,   projected Veronese surfaces  and quadruple planes can be deformed to Veronese surfaces in a 5-dimensional quadric. We give some concrete  examples of this in (\ref{n>4.main.thm.exmp.2}).

For now we focus on 
 describing the birational type of
$\maps_V^{\circ}(\p^2, Q^n)$.

\begin{say}[Proof of Theorem \ref{n>4.main.thm}]\label{n>4.main.pf}
We closely follow the arguments in (\ref{Ver.final.pf.say}).
As in (\ref{Ver.final.pf.say}.4), the key point is to  study  the  correspondence
$$
\{V\}\leftarrow  \{V\cup Q^2\cup P_1\cup P_2\subset Q^4\}\rightarrow \{P_1\cup P_2\}.
\eqno{(\ref{n>4.main.pf}.1)}
$$ 
Here $V$ is a Veronese surface in $Q^n$, $Q^4\subset Q^n$ is a smooth subquadric, $V\cup Q^2\cup P_1\cup P_2 $ is a complete intersection of 3 quadrics in $\p^{n+1}$ and 
 $V\cup Q^2\cup P_1\cup P_2 \subset Q^4$ is the same configuration as in
(\ref{Ver.final.pf.say}).

Since a  Veronese surface spans a $\p^5$, it uniquely determines the
4-dimensional subquadric
$Q^4\subset Q^n$ that contains it.
 Thus, as before,  the left hand side of (\ref{n>4.main.pf}.1)
is birational to
$$
\chow_V(Q^n)\leftarrow  \chow_V(Q^n) \times \p^4.
\eqno{(\ref{n>4.main.pf}.2)}
$$ 
On the right hand  side of (\ref{n>4.main.pf}.1)  we start with 
the family of pairs of intersecting planes 
 $\{P_1\cup P_2\} $. Let us denote this by 
$\sym^2_i\bigl( \OG(\p^2, Q^n)\bigr)$.
(This is a locally closed subset of 
$\sym^{2}\bigl( \OG(\p^2, Q^{4})\bigr)$, but,
unlike in the $n=4$ case, it is not an irreducible component.
However, it is a dense subset of  an irreducible component
of   $\hilb_{t^2+3t+1}(Q^n)$, parametrizing subschemes with Hilbert polynomial $t^2+3t+1 $.)

Sending a pair $P_1\cup P_2$ to the intersection point gives a map
$$
\Pi:\sym^2_i\bigl( \OG(\p^2, Q^n)\bigr)\to Q^n
\eqno{(\ref{n>4.main.pf}.3)}
$$
whose fiber over $p\in Q^n$ is a dense open subset of
$\sym^2\bigl( \OG(\p^1, Q^{n-2}_p)\bigr)$.
In general the fibration given by $\Pi$ is not birationally trivial, but
it is so if $Q^n(k)\neq\emptyset$; this follows from   (\ref{basic.not.not}.3). We obtain that if $Q^n(k)\neq\emptyset$ then
$$
\sym^2_i\bigl( \OG(\p^2, Q^n)\bigr)\simb Q^n\times \sym^2\bigl( \OG(\p^1, Q^{n-2}_W)\bigr).
\eqno{(\ref{n>4.main.pf}.4)}
$$
The rest is very much like before, except that
the linear span of $P_1\cup P_2$ is  4-dimensional, thus the choice of
$Q^4\supset P_1\cup P_2$ is equivalent to choosing a 5-dimensional
linear subspace containing $P_1\cup P_2$.  
So we get an extra factor of $ \p^{n-4}$. 
At the end  we obtain that
 the diagram (\ref{n>4.main.pf}.1)
is birational to
$$
\begin{array}{ccc}
\chow_V(Q^n) \times \p^4
&\simb& \sym^{2}_i\bigl( \OG(\p^2, Q^{n})\bigr) \times \p^{n-4}\times \p^{10}\\
\downarrow && \downarrow\\
\chow_V(Q^n) &&
 \sym^{2}_i\bigl( \OG(\p^2, Q^{n})\bigr).
\end{array}
\eqno{(\ref{n>4.main.pf}.5)}
$$ 
If $Q^n(k)\neq\emptyset$ then (\ref{n>4.main.pf}.4)
further gives that
$$
\chow_V(Q^n)\times \p^4\simb 
\sym^2\bigl( \OG(\p^1, Q^{n-2}_W)\bigr)\times \p^{2n+6}.
\eqno{(\ref{n>4.main.pf}.6)}
$$
Since $\maps_V^{\circ}(\p^2,Q^n) \simb \chow_V(Q^4)\times \p^8$, we see that
(\ref{n>4.main.pf}.6) implies (\ref{n>4.main.thm}.2).
As we already noted,  (\ref{pls.say}.2)  shows that 
(\ref{n>4.main.thm}.1) and (\ref{n>4.main.thm}.2) are equivalent.
\qed

\end{say}

\begin{exmp} \label{n>4.main.thm.exmp.2}
We give examples
of   projected Veronese surfaces  and quadruple planes  deforming to Veronese surfaces in a 5-dimensional quadric. 

First consider  the quadric $Q^5:=
\bigl(x_0x_1+x_2x_3=x_4^2+ax_5^2+x_6^2\bigr)$.
Pick any degree 2 polynomial $q_1(u,v)$ that is relatively prime to 
$u^2+av^2$.  Pick any degree 2 polynomial
$q_6(u,v)$ and write $q_6^2=q_0q+q_3(u^2+av^2)$ where $q_0, q_3$ have degree 2. 
Consider the family of
maps 
$$
\phi_t:(u{:}v{:}w)\mapsto \bigl(t^2q_0: q_1: u^2+av^2: w^2+t^2q_3: uw:vw:tq_6\bigr)
$$
Note that  $\phi_0$ is a  projected Veronese surface but the image of $\phi_t$ is a Veronese surface for $t\neq 0$.

Second, for the quadric  $Q^5:=(x_0x_3+x_1x_4+x_2x_5=x_6^2)$ consider the family of
maps $\phi_t:\p^2\to  Q^5$ given by
$$
\phi_t:(u{:}v{:}w)\mapsto \bigl(u^2 : v^2: w^2:(tv+tw)^2:(tu+tw)^2:(tu+tv)^2:
tuv+tuw+tvw\bigr).
$$
Note that  $\phi_0$ is a quadruple plane but the image of $\phi_t$ is a Veronese surface for $t\neq 0$.
\end{exmp}
 
\section{Orthogonal Grassmannians}\label{sec.OG}

In the next two sections   we discuss various facts about
  orthogonal Grassmannians that we used earlier.
For the current applications we  mainly need symmetric squares of 
$\OG(\p^r, Q^n) $ for $r=1,2$. We focus on these cases but give
more general statements when possible.

\begin{defn}[Orthogonal Grassmannians] \label{OG.defn}
The scheme parametrizing  of all $m$-dimensional linear spaces contained in a smooth quadric   $Q^n$ is called 
the {\it orthogonal Grassmannian,} denoted by $\OG(\p^m,Q^n)$.

Over an algebraically closed field the orthogonal Grassmannian
is usually denoted by 
$\OG(m+1,n+2)$;
the notation suggests $(m+1)$-dimensional $q$-isotropic sub-vector spaces
of  $k^{n+2}$ where $q$ is a non-degenerate quadratic  form.

The orthogonal Grassmannian $\OG(m+1,n+2)$ is a smooth projective
variety that is homogeneous under the  orthogonal group $\OO(q)$.

The dimension of  orthogonal Grassmannians can be computed several ways.
For example, an easy application of (\ref{OG.say}.2)
shows that the universal $\p^m$-bundle over $\OG(m+1,n+2) $
is an $\OG(m,n) $-bundle over $Q^n$. 
Thus 
$\dim \OG(m+1,n+2)=(n-m)+\dim \OG(m, n)$.
Repeating this until we reach  $\OG(1,n-2m+2)\cong Q^{n-2m}$ we get that 
$$
\dim \OG(m+1,n+2)=\tfrac12 (2n-3m)(m+1).
\eqno{(\ref{OG.defn}.1)}
$$
This also shows that $\OG(m+1,n+2)$
 is nonempty iff $2m\leq n$,  irreducible for $2m<n$ and 
 has 2 irreducible  components for $2m=n$; see (\ref{mid.OG.comps}).
(Over an algebraically closed field some authors use  $\OG(m+1,2m+2)$ to denote any one of these components.)

Thus $\OG(\p^m,Q^n)$ is a smooth projective
$k$-variety that is homogeneous under the group $\PGO(Q^n)=\aut(Q^n)$.
\end{defn}

Intersection theory of  orthogonal Grassmannians is studied in 
\cite{MR3219572} and their quantum cohomology in \cite{MR2027200, MR3420327}. 
Here we aim to study  the  birational properties of orthogonal Grassmannians
and their symmetric powers. 
Let us start with the existence of $k$-points.

\begin{say}[Basic correspondences] \label{OG.say}

A basic tool comparing different orthogonal Grassmannians is the variety parametrizing pairs of linear spaces
$L^r\subset L^m$ contained in $ Q^n$.
Let us denote it by $\OF(\p^r, \p^m, Q^n)$ (for orthogonal flags). 
Forgetting one of the linear spaces gives two morphisms
$$
 \OG(\p^r,Q^n)\stackrel{\pi_r}{\longleftarrow}
  \OF(\p^r, \p^m, Q^n)\stackrel{\pi_m}{\longrightarrow}   \OG(\p^m,Q^n).
\eqno{(\ref{OG.say}.1)}
$$
The fiber of $\pi_m$ over  $L^m\subset Q^n$ is the (ordinary)
Grassmannian  $\grass(\p^r, L^m)$.

We claim that  the 
fiber of $\pi_r$ over  $L:=L^r\subset Q^n$ is an
orthogonal Grassmannian $\OG(\p^{m-r-1}, Q^{n-2r-2}_L)$ that can be constructed as follows. First choose coordinates such that 
$$
q(x_0,\dots, x_{n+1})\sim y_0y_1+\cdots +y_{2r}y_{2r+1}+
q_{n-2r}(y_{2r+2},\dots, y_{n+1})
$$
and $L^r=(y_1=y_3=\cdots=y_{2r+1}=y_{2r+2}=y_{2r+3}=\cdots = y_{n+1}=0)$. 
Then we can identify an $L^m$ satisfying $L^r\subset L^m\subset Q^n$
with a linear space
$$
L^{m-r-1}\subset Q^{n-2r-2}_{L}:=\bigl(q_{n-2r}(y_{2r+2},\dots, y_{n+1})=0\bigr).
$$
Geometrically one can obtain $ Q^{n-2r-2}_{L} $ as follows.
First take
$$
 Q^{n-r-1,r}_{L}:=\cap_{p\in L}T_pQ^n.
$$
It is a quadric of dimension $n-r-1$ that is singular along $L$.
Then $Q^{n-r-1,r}_{L}$ is the cone over $ Q^{n-2r-2}_{L} $ with vertex $L$.

Two   cases are especially useful. If $r=0$ then we get
$$
Q^n\stackrel{\pi_0}{\longleftarrow}
  \OF(\p^0, \p^m, Q^n)\stackrel{\pi_m}{\longrightarrow}   \OG(\p^m,Q^n),
\eqno{(\ref{OG.say}.2)}
$$
where the fiber of $\pi_m$ over  $L^m\subset Q^n$ is $ L^m$ while the 
fiber of $\pi_0$ over  $p\in Q^n$ is the
orthogonal Grassmannian $\OG(\p^{m-1}, Q^{n-2}_p)$. 

The other interesting case is $r=m-1$. Then we have
$$
 \OG(\p^{m-1},Q^n)\stackrel{\pi_{m-1}}{\longleftarrow}
  \OF(\p^{m-1}, \p^m, Q^n)\stackrel{\pi_m}{\longrightarrow}   \OG(\p^m,Q^n).
\eqno{(\ref{OG.say}.3)}
$$
Here  $\pi_m$ is a 
$\p^m$-bundle  (the dual of the natural $\p^m$-bundle)
 while the 
fiber of $\pi_{m-1}$ over  $L^{m-1}\subset Q^n$ is the
quadric $Q^{n-2m}_L)$.
\medskip

{\it Claim \ref{OG.say}.4.} Using the notation of 
 (\ref{OG.say}.1),
\begin{enumerate}
\item $\pi_m$  is a Zariski locally trivial $\p^m$-bundle and
\item $\pi_r$  is a Zariski locally trivial $\OG(\p^{m-r-1}, Q^{n-2r}_M)$-bundle
if there is a linear subspace  $M^r\subset Q^n$ defined over $k$.
\end{enumerate}
\medskip

Proof.  
If $r=0$ then  $\pi_m$ is a $\p^m$-bundle and the pull-back of
$\o_{\p^{n+1}}(1)$ gives a relative $\o(1)$. Thus the $\p^m$-bundle is
Zariski locally trivial and so are its Grassmann bundles, proving the first part.

By contrast, $\pi_r$ is usually not  Zariski locally trivial.
However,  assume that we have   $M:=M^r\subset Q^n$ defined over $k$.
As in  (\ref{basic.not.not}.2), for any $L:=L^r\subset  Q^n\setminus Q^{n-r-1,r}_{M}$ 
we get natural isomorphisms  (defined over $k(L)$)
$$
Q^{n-2r-2}_L\cong Q^{n-r-1,r}_L\cap Q^{n-r-1,r}_{M}\cong Q^{n-2r-2}_{M}.
$$
These trivialize the $Q^{n-2r-2}_L$-bundle over the open set 
parametrizing those $L$ that are disjoint from $Q^{n-r-1,r}_{M}$. 
Thus the $\OG(\p^{m-r-1}, Q^{n-2r-2}_L)$-bundle is also trivialized.
\qed

\end{say}

\begin{prop}\label{bir.to.ani.prop.prop}
 Let $Q^n$ be a smooth quadric such that $Q^n(k)\neq \emptyset$. Then  
$$
\OG(\p^m,Q^n)\simb \OG(\p^{m-1}, Q^{n-2}_W)\times \p^{n-m}.
$$
\end{prop}

Proof. Pick a $k$ point $p\in Q^n$ and 
consider the diagram (\ref{OG.say}.2)
$$
Q^n\stackrel{\pi_0}{\longleftarrow}
  \OF(\p^0, \p^m, Q^n)\stackrel{\pi_m}{\longrightarrow}   \OG(\p^m,Q^n).
$$
Let $Q^{n-m}\subset Q^n$ be the intersection of $Q^n$ with a 
linear subspace of codimension $m$ containing $p$ such that $Q^{n-m}$ is smooth at $p$.
Since the 
 fiber of $\pi_m$ over  $L^m\subset Q^n$ is $ L^m$,
we see that the restriction
$$
\pi_m: \pi_0^{-1}(Q^{n-m})\to  \OG(\p^m,Q^n)
$$
is birational. By (\ref{OG.say}.4), $\pi_0$ is Zariski locally trivial with
fiber  $\OG(\p^{m-1}, Q^{n-2}_W)$.  Projecting from $p$ shows that  $Q^{n-m}\simb  \p^{n-m}$. Thus
$$
\OG(\p^m,Q^n)\simb \OG(\p^{m-1}, Q^{n-2}_W)\times Q^{n-m}
\simb \OG(\p^{m-1}, Q^{n-2}_W)\times \p^{n-m}. \qed
$$

\begin{cor}\label{bir.to.ani.prop}
 Let $Q^n$ be a quadric with Witt index $r$  and
anisotropic reduction $Q^{n-2r}_{\rm an}$.  
\begin{enumerate}
\item If $m< r$ then $\OG(\p^m, Q^n)$ is rational.
\item If $m\geq r$ then $\OG(\p^m, Q^n)(k)=\emptyset$  and 
$$
\OG(\p^m, Q^n)\simb \OG(\p^{m-r}, Q^{n-2r}_{\rm an})\times \p^{r(2n-r+1)/2}. 
$$
\end{enumerate}
\end{cor}

Proof. The claims are clear  if the Witt index is $0$. 
Otherwise we apply (\ref{bir.to.ani.prop.prop}) repeatedly.

If $m\geq  r$ we stop when we reach the anisotropic kernel
$Q^{n-2r}_{\rm an}$ to get (2). If $m<r$ then we stop with
$\OG(\p^{0}, Q^{n-2m}_W)=Q^{n-2m}_W$. Since $m<r$ the quadric
$Q^{n-2m}_W$ still has $k$-points, hence it is birational to $\p^{n-2m}$.
Thus  $\OG(\p^m, Q^n)$ is rational. \qed
\medskip

We get the following rationality criteria for
$\OG(\p^m,Q^n)$.

\begin{thm} \label{pt.rtl.cor}
Let $Q^n=\bigl(q(x_0,\dots, x_{n+1})=0\bigr)$  be a smooth quadric
over a field $k$. 
The following are equivalent.
\begin{enumerate}
\item  $\OG(\p^m,Q^n)$ is rational.
\item  $\OG(\p^m,Q^n)$ has a $k$-point.
\item  $\OG(\p^m,Q^n)$ has a $k'$-point  over
an odd degree field extension  $k'/k$.
\item $Q^n$ contains a linear subspace of dimension $m$. 
\item The Witt index of $q$ is $\geq m+1$.
Equivalently, 
we can write
$$
q(x_0,\dots, x_{n+1})\sim y_0y_1+\cdots +y_{2m}y_{2m+1}+
q_W(y_{2m+2},\dots, y_{n+1}). 
$$
\end{enumerate}
\end{thm}

Proof. The implications  (1) $\Rightarrow$ (2) $\Rightarrow$ (3)
and (2) $\Leftrightarrow$ (4) are clear and
(4) $\Leftrightarrow$ (5) was noted in Paragraph \ref{basic.not.not}. 
We proved in (\ref{bir.to.ani.prop}) that (5) $\Rightarrow$ (1).
Thus it remains to show that (3) $\Rightarrow$ (2).

 We use induction on $n$. If  $\OG(\p^m, Q^n) $  has a $k'$-point  then $Q^n$ contains a linear space defined over $k'$,
hence also a $k'$-point. Thus there is a $k$-point $p\in Q^n$ by
Springer's theorem; see \cite{MR0047021} or (\ref{sym.Q.thm}). By (\ref{bir.to.ani.prop.prop})
$$
\OG(\p^m, Q^n)\simb \OG(\p^{m-1}, Q^{n-2}_W)\times \p^{n},
\eqno{(\ref{pt.rtl.cor}.6)}
$$
hence $\OG(\p^{m-1}, Q^{n-2}_W) $  also has a $k'$-point. By induction
 $\OG(\p^{m-1}, Q^{n-2}_W) $  has a $k$-point, hence,
using (\ref{pt.rtl.cor}.6) from right to left shows that
 $\OG(\p^{m}, Q^{n}) $  has a $k$-point. \qed
\medskip

 Three of the  orthogonal Grassmannians
are closely related to Severi-Brauer varieties.
The obvious one is $\OG(\p^1, Q^2) $, the more interesting
examples are $\OG(\p^1, Q^3) $ and $\OG(\p^2, Q^4) $.
For these the
symmetric powers  can be determined using \cite{k-spsb}.
(See \cite{k-sbnotes} for a geometric introduction to Severi-Brauer varieties.)

\begin{lem} \label{og13.sb.thm}
Let  $Q^3=\bigl(q(x_0,\dots, x_4)=0\bigr)$ be  a smooth quadric 3-fold over $k$. Then
\begin{enumerate}
\item $\OG(\p^1, Q^3) $ is a 3-dimensional Severi--Brauer variety.
\item $\OG(\p^1, Q^3) $ has period 1 or 2.
\item $\OG(\p^1, Q^3)\cong \p^3 $ $\Leftrightarrow$   $Q^3$ contains a line   $\Leftrightarrow$ $q\sim y_0y_1+y_2y_3+a_4y_4^2$.
\item  $\OG(\p^1, Q^3) $ has index 2 iff   $q\sim q_4(y_0,\dots, y_3)+a_4y_4^2$
such that $\disc(q_4)$ is a square. If these hold then 
 $\sym^2\bigl(\OG(\p^1, Q^3)\bigr)\simb \p^6$.
\item  $\OG(\p^1, Q^3) $ has index 4 iff
 $\sym^2\bigl(\OG(\p^1, Q^3)\bigr)$ has no $k$-points.
\item Let $a_0x_0^2+\cdots+a_4x_4^2=0$ be a diagonal equation of $Q^3$. Then
$$
\bigl[\OG(\p^1, Q^3)\bigr]=\bigl[(a_0x_0^2+a_1x_1^2+a_2x_2^2=0)\bigr]
\times \bigl[(a_0a_1a_2y^2=a_3x_3^2+a_4x_4^2)\bigr],
$$
where $\bigl[*\bigr]$ denotes the class of an object in the Brauer group.
\end{enumerate}
\end{lem}

Proof. Assume first that $k$ is algebraically closed and
 let $\ell\in Q^3$ be any line. Then
$\{\ell':\ell\cap \ell'\neq\emptyset\}$ defines a divisor $H_{\ell}\subset \OG(2,5)$. Given three general lines $\ell_1, \ell_2, \ell_3$, there is a unique line $\ell'$ that meets all 3. Indeed, the first 2 lines span an
$L^3\subset \p^4$ and  $\ell_3\cap L^3$ is a single point $p$. Thus the only line $\ell'$ in $\p^4$  that meets the $\ell_i$ is the one passing through
$p$ and meeting both  
$\ell_1, \ell_2$. This $\ell'$ meets $Q^3$ in 3 points,so is contained in it.
This shows that  $(H^3)=1$, so  $|H|$ maps $ \OG(2,5)$
isomorphically onto $\p^3$, proving (1). 

An $m$-dimensional Severi-Brauer variety is isomorphic to  $\p^m$
iff it has a $k$-point. Thus $\OG(\p^1, Q^3)\cong \p^3 $ iff  $Q^3$ contains a line and the rest of (3) is clear.

If we are over a field $k$, then there may not be lines defined over $k$ but there are conics $C\subset Q^3$ defined over $k$.
Thus  
$\{\ell':C\cap \ell'\neq\emptyset\}$ defines a divisor $D_{C}\subset \OG(\p^1, Q^3)$ that is defined over $k$. Over $\bar k$ we have $D_C\sim 2H$
and $\ell'\mapsto C\cap \ell'$ gives a morphism $D_C\to C$ defined over $k$.
Thus in fact $D$ is a product of 2 conics $D_C\cong C\times C'$.
This shows that  $\OG(\p^1, Q^3)$ is a  Severi-Brauer variety 
of period 2, proving (2). (By our construction, it contains a quadric that is a product of 2 conics; actually this always holds, by a theorem of Albert, if the period is 2; see
\cite[Thm.18]{k-conic} for a geometric proof.)

By definition, $\OG(\p^1, Q^3)$ has index 2 iff it  has a degree 2 point,
this $Q^3$ contains a conjugate pair of lines $\ell, \ell'$.
These span an $L^3\subset \p^4$ and $Q^2:=L^3\cap Q^3$ is a
smooth quadric that contains a  conjugate pair of lines.
By (\ref{mid.OG.comps}) this holds iff its discriminant is a square.

Next we claim that if $P$ is a Severi-Brauer variety that has a  point 
over a degree 2 field extension $k'/k$ then $\sym^2P$ is rational. Indeed, let $\pi:U\to G$
denote the universal bundle over the Grassmannian of lines. 
Then $\sym^2P$ is birational to  $\p_G\bigl(\pi_*\o_U(-K_U)\bigr)$, hence to
 $\p^2\times G$. Let $L$ be a hyperplane in $P$ defined over $k'$. Then $G$ is birational to the Weil restriction $\Re_{k'/k}L$, hence rational and so is
$\sym^2 P$.  This proves (4) and (5) just lists the remaining cases.

For (6) we need to follow the construction of $D_C$ above. So assume that  $Q^3=(a_0x_0+\cdots+a_4x_4^2=0)$. We can use  
$C=(x_3=x_4=0)$ as our first conic; it is isomorphic to
$(a_0x_0^2+a_1x_1^2+a_2x_2^2=0)\subset \p^2$.
To get our second conic $C'$, we project $Q^3$ from $C$ to get
$\pi:Q^3\map \p^1_{x_3,x_4}$. The  fiber over a point $(s{:}t)$ is 
$$
\begin{array}{rcl}
Q^2(s,t)&=&\bigl(a_0x_0^2+\cdots+a_4x_4^2=tx_3-sx_4=0\bigr)\\
&\cong& \bigl(a_0x_0^2+a_1x_1^2+a_2x_2^2+(a_3s^2+a_4t^2)x_3^2=0\bigr).
\end{array}
$$ 
Furthermore,  $C'\to \p^1_{x_3,x_4}$ is the double cover whose fiber 
over $(s{:}t)$ is the 2 families of lines on $Q^2(s,t)$.
Thus, by (\ref{mid.OG.comps}),  the equation of $C'$ is  $w^2=\disc\bigl( Q^2(s,t)\bigr)$.
We can divide  by $ a_0a_1a_2$ and set $y:=w/(a_0a_1a_2)$ to get the form
$C'\cong (a_0a_1a_2y^2=a_3x_3^2+a_4x_4^2)$. \qed

\medskip

As we noted, the middle dimensional  orthogonal Grassmannians $\OG(\p^{r}, Q^{2r})$ are geometrically reducible. The discriminant tells us when
they are reducible.

\begin{lem} \label{mid.OG.comps}
Let $Q^{2r}$ be an even dimensional smooth quadric over
 a field $k$.
Then the scheme parametrizing the 2 geometric components of
$\OG(\p^{r}, Q^{2r})$ is  $\bigl(t^2+(-1)^r\Delta=0\bigr)$,
where $\Delta$ is the discriminant of $Q^{2r}$. 

Thus $\OG(\p^{r}, Q^{2r})$ is reducible iff
$(-1)^{r+1}\Delta(Q^{2r})$ is a square.
\end{lem}

Proof. Set $m:=r+1$ and consider the universal quadratic form
$$
q^{2m}_{\rm univ}(x_1,\dots, x_{2m}):=\sum_{1\leq i\leq j\leq 2m} t_{ij}x_ix_j
$$
over the affine space of dimension $m(2m+1)$ over $\spec \z[\tfrac12]$.
Let $\Delta$ be the discriminant and 
$U:=\a^{m(2m+1)}\setminus (\Delta=0)$ the complement of its zero set.
 The Stein factorization of the universal orthogonal Grassmannian
of $m$-dimensional isotropic subspaces gives an \'etale double cover of $U$.
Since $\a^{m(2m+1)} $ is simply connected, this is given by an equation
$$
z^2=(-1)^a\cdot 2^b\cdot  \Delta^c
$$
for some integers  $a,b,c$ which can all be chosen  $0$ or $1$.
To pin down these choices, we can check a particular case, for instance
$$
q_t:=x_1^2-x_2^2+\dots+x_{2m-3}^2-x_{2m-2}^2+x_{2m-1}^2-tx_{2m}^2.
$$
We have $\Delta(q_t)=(-1)^mt$ and the 2 families of $m$-planes
correspond to the members
$$
\bigl(x_1-x_2=\dots=x_{2n-3}-x_{2n-2}=x_{2m-1}\pm\sqrt{t}x_{2m}=0.
$$
Thus the correct equation is $z^2=(-1)^m\Delta$. \qed

\begin{prop} \label{mid.OG.comps.red}
Let $Q^{2r}$ be an even dimensional smooth quadric over
 a field $k$ and $Q^{2r-1}\subset Q^{2r}$ a smooth quadric hypersurface.
Them
$$
\OG(\p^{r}, Q^{2r})\cong \OG(\p^{r-1}, Q^{2r-1})\times 
\bigl(t^2+(-1)^r\Delta=0\bigr).
$$
\end{prop}

Proof. Intersecting $L^r\subset Q^{2r}$ with $Q^{2r-1}$ gives a
morphism $ \pi_1: \OG(\p^{r}, Q^{2r})\to \OG(\p^{r-1}, Q^{2r-1})$
which is everywhere defined since a smooth quadric of dimension $2r-1$ can not contain a $\p^r$. We can think of a fiber of $\pi_1$ as a
fiber of the projection considered in (\ref{OG.say}.3)
$$
 \OG(\p^{r-1},Q^{2r})\stackrel{\pi_{r-1}}{\longleftarrow}
  \OF(\p^{r-1}, \p^r, Q^{2r}).
$$
We identified this fiber as  a 0-dimensional quadric, hence 2 points.
Combining with (\ref{mid.OG.comps}) we get a surjection
$$
\OG(\p^{r}, Q^{2r})\to \OG(\p^{r-1}, Q^{2r-1})\times \bigl(t^2+(-1)^r\Delta=0\bigr),
$$
which is then necessarily an isomorphism.\qed

\begin{cor} $\OG(\p^2, Q^4) $ is a Severi-Brauer scheme over
$\bigl(t^2+\Delta=0\bigr) $.\qed
\end{cor}

\section{Symmetric squares of  orthogonal Grassmannians}
\label{sec.ssog}

In this section we  study the rationality of $\sym^{2}\bigl(\OG(\p^1, Q^n)\bigr)$.
Note that if $n\geq 3$ then $\aut(Q^n)$ acts with a dense orbit on
$\sym^{2}\bigl(\OG(\p^1, Q^n)\bigr)$. Thus if the latter has a $k$-point then it is   unirational. (This needs a small argument about singular $k$-points.)

More generally, one can ask about the structure of 
 $\sym^{d}\bigl(\OG(\p^m, Q^n)\bigr)$.
Building on the assumption that quadrics are similar to 
Severi-Brauer varieties,
the following questions appear to be quite natural and they should give a general outline of what could hold. The analogous assertions  for
Severi-Brauer varieties are proved in \cite{k-spsb},  however,    similar  problems
about other homogeneous spaces  turned out to have a negative answer 
\cite{MR2097288, MR2176607}.
So there is little reason to believe that the answers should be positive.

\begin{ques} \label{lin.main.symd.ques} Let $Q^n$ be a smooth quadric  over a field $k$.
Fix $m\geq 0$ and 
write $d=2^ad'$ where $d'$ is odd. Which of the following are true?
\begin{enumerate}
\item $\sym^{d}\bigl(\OG(\p^m, Q^n)\bigr)$ is rational over $k$ if $2^{m+1}\mid d$.
\item  $\sym^{d}\bigl(\OG(\p^m, Q^n)\bigr)$ is rational over $k$ iff it has a $k$-point.
\item The stable birational equivalence class of $\sym^{d}\bigl(\OG(\p^m, Q^n)\bigr)$ depends only on $2^a$. 
\end{enumerate}
\end{ques}

The simplest case, when $a=0$, follows from (\ref{pt.rtl.cor}).

\begin{prop}  \label{lin.main.symd.odd} Let $Q^n$ be a smooth quadric  over a field $k$ and assume that $d$ is odd. Then
$$
\sym^{d}\bigl(\OG(\p^m, Q^n)\bigr)\sims \OG(\p^m, Q^n).
$$
\end{prop}

Proof. We show the stable birational equivalences
$$
\sym^{d}\bigl(\OG(\p^m, Q^n)\bigr)\sims 
\sym^{d}\bigl(\OG(\p^m, Q^n)\bigr)\times \OG(\p^m, Q^n)
\sims \OG(\p^m, Q^n).
$$
First let $K$ be the function field of $ \OG(\p^m, Q^n)$.
Then $ \OG(\p^m, Q^n)(K)\neq \emptyset$ hence
$ \OG(\p^m, Q^n)$ is rational over $K$ by (\ref{pt.rtl.cor}) and so are its symmetric powers. This shows the stable birational equivalence on the right.

Next let $L$ be the function field of $ \sym^{d}\bigl(\OG(\p^m, Q^n)\bigr)$.
Then $\OG(\p^m, Q^n)$ has a point in an odd degree field extension $L'/L$,
hence it  is rational over $L$  by (\ref{pt.rtl.cor}). Thus the   
stable birational equivalence on the left also holds. \qed

\medskip

Another special case of Question \ref{lin.main.symd.ques} is 
symmetric powers of $\OG(\p^0,Q^n)\cong Q^n$.
\cite{MR0047021} proves that if a quadric has a point over
an odd degree field extension of $k$ then it has a $k$-point.
The idea of the proof gives more.

\begin{thm}\label{sym.Q.thm}
Let $Q^n$ be a smooth quadric of dimension $n\geq 1$. Then 
$$
\sym^d(Q^n)\simb 
\left\{
\begin{array}{ll}
\p^{nd}&\qtq{if $d$ is even,}\\
Q^n\times \p^{nd-n}&\qtq{if $d$ is odd.}
\end{array}
\right.
$$
\end{thm}

Proof. We prove that, for $d\geq 2$,
$\sym^d(Q^n)\simb \sym^{d-2}(Q^n)\times \p^{2n}$.
In affine coordinates we write
$Q^n=\bigl(q(x_1,\dots,x_{n},t)=0\bigr)$.   

For $d$ general points $p_1,\dots, p_d\in Q^n$ there are unique polynomials 
$g_1(t),\dots,g_n(t)$ 
of degree $\leq d-1$ whose graph $\Gamma:=\Gamma(p_1,\dots, p_d)$ 
passes through  $p_1,\dots, p_d$.
We define $\sym^d(Q^n)\map \sym^{d-2}(Q^n)\times \a^{n}\times \a^{n}$
by sending  $\{p_1,\dots, p_d\}$ to the triple
consisting of the $d-2$ other intersection points of  $\Gamma$ 
with $Q^n$ and the  intersection points of  $\Gamma$ 
with the hyperplanes $(t=0)$ and $(t=1)$. \qed

\begin{say}[Witt reduction and symmetric powers]\label{pls.say}
The rationality problem of symmetric powers can be reduced to the anisotropic case.
Indeed, we proved in (\ref{bir.to.ani.prop.prop}) that 
if  $Q^n$ is a smooth quadric such that 
 $Q^n(k)\neq \emptyset$ then 
$$
\OG(\p^m, Q^n)
\simb \OG(\p^{m-1}, Q^{n-2}_W)\times  \p^{n-m}. 
\eqno{(\ref{pls.say}.1)}
$$
We can  next take symmetric powers and, using 
\cite[Lem.7]{k-spsb},  conclude that
$$
\sym^d\bigl(\OG(\p^m, Q^n)\bigr)\simb \sym^d\bigl(\OG(\p^{m-1}, Q^{n-2}_W)\bigr)\times  \p^{d(n-m)}. 
\eqno{(\ref{pls.say}.2)}
$$
\end{say}

\begin{say}[Symmetric square of $\OG(\p^{r}, Q^{2r})$]
\label{sq.og.some.say}
Since
$\OG(\p^r, Q^{2r})$ has 2 geometric irreducible components, 
its symmetric square can be written as
$$
\sym^2\bigl( \OG(\p^r, Q^{2r})\bigr)=\sym^{2}_{s}\bigl( \OG(\p^r, Q^{2r})\bigr)\amalg
\sym^{2}_{a}\bigl( \OG(\p^r, Q^{2r})\bigr)
\eqno{(\ref{sq.og.some.say}.1)}
$$
where the first (symmetric) component parametrizes pairs
$L_1^r, L_2^r\subset Q^{2r}$ that are in different geometric irreducible components
of $\OG(\p^r, Q^{2r})$ while the second (anti-symmetric) component  parametrizes pairs
that are in the same geometric irreducible component.
Note that $\sym^{2}_{s}\bigl( \OG(\p^r, Q^{2r})\bigr) $ is
geometrically irreducible while 
$\sym^{2}_{a}\bigl( \OG(\p^r, Q^{2r})\bigr) $ is
geometrically reducible. 

If $r=1$ then 
$$
\begin{array}{lcl}
\sym^{2}_{s}\bigl( \OG(\p^1, Q^{2})\bigr)&\cong &  Q^2 \qtq{and}\\
\sym^{2}_{a}\bigl( \OG(\p^1, Q^{2})\bigr) &\cong& (t^2=\Delta)\times \p^2
\end{array}
\eqno{(\ref{sq.og.some.say}.2)}
$$
Indeed, the first isomorphism is obtained by identifying a
pair of intersecting lines on $Q^2$ with their intersection point.
The second isomorphism is obtained by noting that
$  \OG(\p^1, Q^{2})$ is a conic bundle over $(t^2=\Delta)$
by (\ref{mid.OG.comps})
and the symmetric square of any  conic is $\p^2$ by (\ref{sym.Q.thm}). 

If $r=2$ then 
two general planes  of the same family meet at a single point,
thus we have a map
$$
\sym^{2}_{a}\bigl( \OG(\p^2, Q^4)\bigr)\map Q^4.
\eqno{(\ref{sq.og.some.say}.3)}
$$
The fiber over $p\in Q^4$ can be identified with
$\sym^{2}_{a}\bigl(\OG(\p^1, Q^2_p)\bigr)$. 
The $\OG(\p^1, Q^2_p)$-bundle over $Q^4$ is a conic bundle over an
\'etale double cover of $Q^4$. Since $Q^4$ is simply connected,
the \'etale double cover is isomorphic to
$Q^4\times (t^2+\Delta=0)$. 
(Note that Witt reduction from $Q^4$ to $Q^2_p$ changes the sign of the discriminant.) Thus we have a conic bundle over
$Q^4\times (t^2+\Delta=0)$ and we take its relative symmetric square.
This is a Zariski locally trivial $\p^2$-bundle over $Q^4\times (t^2+\Delta=0)$, hence
$$
\sym^{2}_{a}\bigl( \OG(\p^2, Q^4)\bigr)\simb  
(t^2+\Delta=0)\times Q^4\times \p^2.
\eqno{(\ref{sq.og.some.say}.4)}
$$
\end{say}

Symmetric squares of $\OG(\p^1, Q^n) $ are related to double covers of Grassmannians. 

\begin{prop}\label{sym2.OG.Qn.cor}
 Let $Q^n$ be a smooth quadric of dimension $\geq 3$. Then there is a double cover
$\pi: G(Q^n)\to \grass(\p^3,  \p^{n+1})$, ramified along a divisor 
$D\sim 2H$
such that 
$$
\sym^2\bigl(\OG(\p^1, Q^n)\bigr)\simb G(Q^n) \times \p^2.
$$
\end{prop}

Proof.
A $K$-point of  $\sym^2\bigl(\OG(\p^1, Q^n)\bigr) $ is the same as a
pair of lines  $\ell,\ell'\subset Q^n$ that are conjugate over $K$.
In general  $\ell,\ell'$ are disjoint, so their span  $\langle \ell, \ell'\rangle$ has dimension 3. Thus we get a map
$$
\pi: \sym^2\bigl(\OG(\p^1, Q^3)\bigr) \map \grass(\p^3,  \p^{n+1}).
\eqno{(\ref{sym2.OG.Qn.cor}.1)}
$$
Given $L^3\subset  \p^{n+1}$ such that $Q^n\cap L^3$ is a smooth quadric surface $Q^2_L$, we proved in (\ref{sq.og.some.say}.2) that
the fiber of $\pi$ over $[L^3]$ is
$$
\sym^{2}_{a}\bigl( \OG(\p^1, Q^{2})\bigr) \cong
 (t^2=\Delta)\times \p^2.
$$
Therefore the Stein factorization of $\pi$ is
$$
\sym^2\bigl(\OG(\p^1, Q^3) \simb  G(Q^n)\times \p^2\to G(Q^n)\to \grass(\p^3,  \p^{n+1}),
\eqno{(\ref{sym2.OG.Qn.cor}.2)}
$$
where $G(Q^n)\to \grass(\p^3,  \p^{n+1}) $ is a double cover
that ramifies over a point of the Grassmannian $[L^3]\in \grass(\p^3,  \p^{n+1})$
iff $Q^n\cap L^3$ is a singular quadric (or if $L^3\subset Q^n$).

Let $Q^3\subset \p^4$ be a smooth quadric 3-fold.  A general pencil of hyperplane sections
$\langle L^3_{\lambda}\cap Q^3\rangle$ of $Q^3$  has 2 singular members.
Since such pencils  $\langle L^3_{\lambda}\rangle\subset L^4\subset  \p^{n+1}$
correspond to lines on $\grass(\p^3,  \p^{n+1}) $, we conclude that
the branch divisor
$D\subset \grass(\p^3,  \p^{n+1})$ is linearly equivalent to
$2H$ where $H$ is the hyperplane class on $\grass(\p^3,  \p^{n+1})$.
\qed

\medskip

Next we show that
$G(Q^3)$ is itself a quadric and  determine its equation. 
This will prove Proposition \ref{sym2.OG.Q3.prop}.

\begin{say}[Proof of Proposition \ref{sym2.OG.Q3.prop}]
\label{sym2.OG.Q3.prop.pf}
By Theorem \ref{n>4.main.thm} and by 
Proposition \ref{sym2.OG.Qn.cor}, 
$$
\maps_2(\p^2, Q^5)\simb \sym^2\bigl(\OG(\p^1, Q^{3}_W)\bigr)\times \p^{20}
\simb G(Q^3_W)\times \p^{22},
$$
where  $Q^3_W$  is the Witt reduction of $Q^5$. Thus we need to
understand $ \pi:G(Q^{3}_W)\to \grass(\p^3,  \p^4)$. 
Note that  $\grass(\p^3,  \p^4)$ is the dual of $\p^4$ and
$\pi$ ramifies along the dual quadric $Q^{\vee}_W$ which is isomorphic to $Q^3_W$.
Thus if $Q^3_W\cong \bigl(q(x_0,\dots, x_4)=0\bigr)$ then 
$$
G(Q^3_W)\cong
\bigl(c\Delta(q)z^2=q(x_0,\dots, x_4)\bigr)\qtq{for some $c$.}
$$
It is not hard to compute everything explicitly and obtain that $c=1$.
As a possible shortcut,
working with the universal quadric over $\z[\frac12]$ as in the proof of
 (\ref{mid.OG.comps}), we see that
$c$ is a unit in $\z[\frac12]$, hence it is enough to compute it in one example. \qed
\end{say}

\begin{ack} I thank I.~Coskun, G.~Di~Cerbo, D.~Krashen, R.~Kusner, M.~Lieblich, N.~Lubbes, R.~Parimala, B.~Poonen,  M.~Skopenkov and B.~Sturmfels
 for
comments, discussions and references. I learned a lot of the early history,
especially the role of  \cite{bertini-book}, from F.~Russo. 
Partial financial support   was provided  by  the NSF under grant number
 DMS-1362960.
\end{ack}


\def\cprime{$'$} \def\cprime{$'$} \def\cprime{$'$} \def\cprime{$'$}
  \def\cprime{$'$} \def\cprime{$'$} \def\cprime{$'$} \def\dbar{\leavevmode\hbox
  to 0pt{\hskip.2ex \accent"16\hss}d} \def\cprime{$'$} \def\cprime{$'$}
  \def\polhk#1{\setbox0=\hbox{#1}{\ooalign{\hidewidth
  \lower1.5ex\hbox{`}\hidewidth\crcr\unhbox0}}} \def\cprime{$'$}
  \def\cprime{$'$} \def\cprime{$'$} \def\cprime{$'$}
  \def\polhk#1{\setbox0=\hbox{#1}{\ooalign{\hidewidth
  \lower1.5ex\hbox{`}\hidewidth\crcr\unhbox0}}} \def\cdprime{$''$}
  \def\cprime{$'$} \def\cprime{$'$} \def\cprime{$'$} \def\cprime{$'$}
\providecommand{\bysame}{\leavevmode\hbox to3em{\hrulefill}\thinspace}
\providecommand{\MR}{\relax\ifhmode\unskip\space\fi MR }
\providecommand{\MRhref}[2]{%
  \href{http://www.ams.org/mathscinet-getitem?mr=#1}{#2}
}
\providecommand{\href}[2]{#2}

\bigskip

\noindent  Princeton University, Princeton NJ 08544-1000

{\begin{verbatim} kollar@math.princeton.edu\end{verbatim}}

\end{document}